% Last Edit 9 September 2002

%   This document is written in Plain TeX
%   The macros: prepictex.tex, pictex.tex, and postpictex.tex are also
%   required for the full compilation of the document.

\magnification=1100
\overfullrule0pt

\input prepictex
\input pictex
\input postpictex
\input amssym.def

% ********************* Definitions ************************************

%\def\widetilde{\mathaccent"0365 }
\def\qed{\hbox{\hskip 1pt\vrule width4pt height 6pt depth1.5pt \hskip 1pt}}

\def\CC{{\Bbb C}}

\def\RR{{\Bbb R}}
\def\ZZ{{\Bbb Z}}

\def\cB{{\cal B}}
\def\cF{{\cal F}}

\def\Hom{{\rm Hom}}
\def\End{{\rm End}}
\def\Ind{{\rm Ind}}

\newsymbol\ltimes 226E % This provides the semidirect product symbol
\newsymbol\rtimes 226F % This provides the other semidirect product symbol

% ********************* FONTS ************************************

\font\smallcaps=cmcsc10
\font\titlefont=cmr10 scaled \magstep1

\font\tinyrm=cmr10 at 8pt
\font\scriptsize=cmr10 at 5pt

% ******************** SECTION HEADERS ***************************

\newcount\sectno
\newcount\subsectno
\newcount\resultno

\def\section #1. #2\par{
\sectno=#1
\resultno=0
\bigskip\noindent\centerline{\smallcaps #1.  #2}~\medbreak}

\def\subsection #1\par{\global\advance\resultno by 1
\medskip\noindent{\bf (\the\sectno.\the\resultno)\ \ #1} }

%******************* MATHEMATICAL LABELS **************************

\def\prop{ \global\advance\resultno by 1
\bigskip\noindent{\bf Proposition \the\sectno.\the\resultno. }\sl}
\def\lemma{ \global\advance\resultno by 1
\bigskip\noindent{\bf Lemma \the\sectno.\the\resultno. }
\sl}
\def\remark{ \global\advance\resultno by 1
\medskip\noindent{\bf Remark \the\sectno.\the\resultno. }}
\def\example{ \global\advance\resultno by 1
\bigskip\noindent{\bf Example \the\sectno.\the\resultno. }}
\def\cor{ \global\advance\resultno by 1
\bigskip\noindent{\bf Corollary \the\sectno.\the\resultno. }\sl}
\def\thm{ \global\advance\resultno by 1
\bigskip\noindent{\bf Theorem \the\sectno.\the\resultno. }\sl}
\def\defn{ \global\advance\resultno by 1
\bigskip\noindent{\it Definition \the\sectno.\the\resultno. }\slrm}
\def\endthm{\rm\bigskip}

\def\endprop{\rm\bigskip}

\def\pf{\rm\bigskip\noindent{\it Proof. }}
\def\endpf{\qed\hfil\bigskip}

\def\Athm{ \global\advance\resultno by 1
\bigskip\noindent{\bf Theorem A.\the\resultno. }\sl}
\def\Alemma{ \global\advance\resultno by 1
\bigskip\noindent{\bf Lemma A.\the\resultno. }
\sl}
\def\Aremark{ \global\advance\resultno by 1
\medskip\noindent{\bf Remark A.\the\resultno. }}

%*************** EQUATIONS WITH NUMBERS **************

\def\Aformula{\global\advance\resultno by 1
\eqno{({\rm A}.\the\resultno)}}
\def\formula{\global\advance\resultno by 1
\eqno{(\the\sectno.\the\resultno)}}
\def\formulano{\global\advance\resultno by 1 (\the\sectno.\the\resultno)}
\def\tableno{\global\advance\resultno by 1
\the\sectno.\the\resultno. }
\def\lformula{\global\advance\resultno by 1
\leqno(\the\sectno.\the\resultno)}

%********** DATING ******************************************
\def\monthname {\ifcase\month\or January\or February\or March\or April\or
May\or June\or
July\or August\or September\or October\or November\or December\fi}

\newcount\mins  \newcount\hours  \hours=\time \mins=\time
\def\now{\divide\hours by60 \multiply\hours by60 \advance\mins by-\hours
      \divide\hours by60         % NOTE: \divide only gives integer answers.
      \ifnum\hours>12 \advance\hours by-12
        \number\hours:\ifnum\mins<10 0\fi\number\mins\ P.M.\else
        \number\hours:\ifnum\mins<10 0\fi\number\mins\ A.M.\fi}

%**************** PAGE HEADERS *************************

\nopagenumbers
\def\runningtitle{\smallcaps affine hecke algebras}
\headline={\ifnum\pageno>1\eoheadline\else\firstheadline\fi}
\def\names{\smallcaps A. Ram\quad and\quad J. Ramagge}
\def\firstheadline{\noindent Preliminary Draft \hfill  \today}
\def\firstheadline{}
\def\eoheadline{\ifodd\pageno\oddheadline\else\evenheadline\fi}
\def\oddheadline{\tenrm\hfil\runningtitle\hfil\folio}
\def\evenheadline{\tenrm \folio\hfil{\names}\hfil}

%**************** TITLE *************************
\vphantom{$ $}  %My kludge to get the first page to move down a bit
\vskip.5truein

\centerline{\titlefont Affine Hecke algebras, cyclotomic Hecke algebras
and Clifford theory}
\bigskip
\centerline{\rm Arun Ram ${}^\ast$}
\centerline{Department of Mathematics}
\centerline{University of Wisconsin, Madison}
\centerline{Madison, WI 53706\ \ USA}
\centerline{{\tt ram@math.wisc.edu}}
%\centerline{Version: \today}
\bigskip
\centerline{\rm Jacqui Ramagge}
\centerline{Department of Mathematics}
\centerline{University of Newcastle}
\centerline{NSW 2308\ \ Australia}
\centerline{{\tt jacqui@maths.newcastle.edu.au}}
%\centerline{Preprint:  April 7, 1999}

\bigskip
\centerline
{\sl Dedicated to Professor C.S. Seshadri on the occasion of
his 70th birthday}

\footnote{}{\tinyrm ${}^\ast$ Research supported in part by NSF
grants DMS-9622985 and DMS-9971099 and the National Security Agency.}

\bigskip

%**************** ABSTRACT *************************
{\it Abstract.}  We show that the Young tableaux theory
and constructions of the irreducible representations
of the Weyl groups of type A, B and D, Iwahori-Hecke algebras of
types A, B, and D, the complex reflection groups $G(r,p,n)$ and
the corresponding cyclotomic Hecke algebras $H_{r,p,n}$,
can be obtained, in all cases, from the affine Hecke algebra
of type A.  The Young tableaux theory was extended to affine Hecke algebras
(of general Lie type) in recent work of A. Ram.
We also show how (in general Lie type) the
representations of general affine Hecke algebras can be constructed from
the representations of simply connected affine Hecke algebras by using
an extended form of Clifford theory.  This extension of Clifford theory
is given in the Appendix.

\section 0. Introduction

Recent work of A. Ram [Ra2,5] gives a
straightforward combinatorial construction
of the simple calibrated modules of affine Hecke algebras (of
general Lie type as well as type A). The first
aim of this paper is to show that Young's seminormal construction and all
of its
previously known generalizations are special cases of the construction
in [Ra5].
In particular, the representation theory of
\smallskip
\item{(a)}  Weyl groups of types A, B, and D,
\smallskip
\item{(b)}  Iwahori-Hecke algebras of types A, B, and D,
\smallskip
\item{(c)}  the complex reflection groups $G(r,p,n)$, and
\smallskip
\item{(d)}  the cyclotomic Hecke algebras $H_{r,p,n}$,
\smallskip\noindent
can be derived
entirely from the representation theory of {\it
affine Hecke algebras of type A.}
Furthermore, the relationship between the affine Hecke algebra
and the objects in (a)-(d) always produces a natural set
of Jucys-Murphy type elements and can be used to prove the
standard Jucys-Murphy type theorems.  In particular, we are able to use
Bernstein's results about the center of the affine Hecke algebra to show
that,
in the semisimple case, the center of the cyclotomic Hecke algebra
$H_{r,1,n}$ is the
set of symmetric polynomials in the Jucys-Murphy elements.

A. Young's seminormal construction of the irreducible representations
of the symmetric group dates from 1931 [Yg1].  Young himself generalized 
his
tableaux to treat the representation theory of
  Weyl groups of types B and D [Yg2].  In 1974 P.N. Hoefsmit [Hf]
generalized the seminormal construction to Iwahori-Hecke algebras of types
A, B, and D.  Hoefsmit's work has never been published and, in 1985,
Dipper and James [DJ, Theorem 4.9] and H. Wenzl [Wz] independently 
treated the seminormal construction for irreducible
representations for Iwahori-Hecke algebras of type A.
In 1994 Ariki and Koike [AK] introduced (some of) the cyclotomic Hecke
algebras and generalized Hoefsmit's construction to these algebras.  The
construction was generalized to a larger class of cyclotomic Hecke
algebras in [Ar2].  For a summary of this work see [Ra1] and [HR].

General affine Hecke algebras are naturally associated to
a reductive algebraic group and the size of the commutative
part of the affine Hecke algebra depends on the structure of the
corresponding algebraic group (simply connected, adjoint, etc.).
The second aim of this paper is to show that it is sufficient to
understand the representation theory of the affine Hecke algebra
in the simply connected case. We describe explicitly how the
representation theory of the other cases is derived
from the simply connected case.

The machine which allows us to accomplish this reduction is a
form of Clifford theory.  Precisely, if $R$ is an algebra
and $G$ is a finite group acting on $R$ by automorphisms
then the representation theory of the ring $R^G$ of
fixed points of the $G$-action can be derived from the
representation theory of $R$ and subgroups of $G$.
This is an extension of the approach to Clifford theory
given by Macdonald [Mac2].

Let us state precisely what is new in this paper.
The main result has three parts:
\smallskip\noindent
\item{(1)} The Hecke algebras $H_{r,p,n}$ of the
complex reflection groups $G(r,p,n)$ can be obtained as fixed point
subalgebras of the Hecke algebra $H_{r,1,n}$ associated to the
complex reflection group $G_{r,1,n}$ via
$$H_{r,p,n} = (H_{r,1,n})^{\ZZ/p\ZZ}.$$
\item{(2)} The Hecke algebras $\tilde H_L$
of nonadjoint $p$-adic groups can be obtained as fixed point
subalgebras of the Hecke algebra $\tilde H_P$
associated to the corresponding adjoint $p$-adic group, via
$$\tilde H_L = (\tilde H_P)^{P/L}.$$
\smallskip\noindent
\item{(3)} There is a form of Clifford theory
(to our knowledge new) that allows one to completely
determine the representation theory of a fixed point subalgebra
$R^G$ in terms of the representation theory of the algebra $R$
and the group $G$.
\smallskip\noindent
We use this method to work out the representation theory of
the algebras $H_{r,p,n}$ in detail.  Additional results include,
\smallskip\noindent
\item{(4)} The discovery of the ``right'' affine
braid groups ${\cal B}_{\infty,p,n}$ and affine Hecke algebras
$H_{\infty,p,n}$ to associate to the cyclotomic Hecke algebras $H_{r,p,n}$.
The representation theory of these new groups and algebras is completely 
determined from the representation theory of the classical
affine braid groups $B_{\infty,1,n}$ and the classical affine
Hecke algebras $H_{\infty,1,n}$ of type $A$ as a consequence of the 
results in (3).
\smallskip\noindent
Finally,
\smallskip\noindent
\item{(5)} We show how the classical trick (due to Cherednik)
for determining the representation theory of the algebras
$H_{r,1,n}$ from that of $H_{\infty,1,n}$ arises from
a map from the affine Hecke algebra of type $C$ to the finite Hecke 
algebra of type $C$ which corresponds to a folding of the Dynkin
diagram.  This explanation is new.  We show that such maps from the 
affine Hecke algebra
to the finite Hecke algebra cannot exist in general type,
and we work out the details of the cases where such homomorphisms do
arise from foldings.

In a recent paper Reeder [Re] has used our results
to prove the Langlands classification
of irreducible representations for general affine
Hecke algebras. (Previously, this was known only in the simply connected 
case, see Kazhdan and Lusztig [KL].)  This also provides a classification 
of the
irreducible constituents of unramified principal series representations
of general split reductive $p$-adic groups. (The Kazhdan-Lusztig
result provides this classification for groups with connected center.)

\bigskip\noindent
{\bf Acknowledgements.}
A. Ram thanks P. Deligne for his encouragement, stimulating questions
and helpful comments on the results in this paper.  We thank D. Passman
for providing some useful references about Clifford theory.
We are grateful for the generous support of this research by the National
Science Foundation, the National Security Agency and
the Australian Research Council.

\section 1.  Algebras with Young tableaux theories

A. Young invented the theory of standard Young tableaux in order to
describe the representation theory of the symmetric group $S_n$; the
group of $n\times n$ matrices such that
\smallskip
\itemitem{(a)}  the entries are either $0$ or $1$,
\smallskip
\itemitem{(b)}  there is exactly one nonzero entry in each row and each 
column.
\smallskip\noindent
Young himself began to generalize the theory and in [Yg1] he provided
a theory for the Weyl groups of type B, i.e. the
hyperoctahedral groups $WB_n\cong (\ZZ/2\ZZ)\wr S_n$ of $n\times n$
matrices such
that
\smallskip
\itemitem{(a)}  the entries are either $0$ or $\pm1$,
\smallskip
\itemitem{(b)}  there is exactly one nonzero entry in each row and each 
column.
\smallskip\noindent
In the same paper Young also treated
the Weyl group $WD_n$ of $n\times n$ matrices such that
\smallskip
\itemitem{(a)}  the entries are either $0$ or $\pm1$,
\smallskip
\itemitem{(b)}  there is exactly one nonzero entry in each row and each 
column,
\smallskip
\itemitem{(c)}  the product of the nonzero entries is $1$.
\smallskip\noindent
W. Specht [Sp] generalized the theory to cover the complex
reflection groups $G(r,1,n)\cong (\ZZ/r\ZZ)\wr S_n$ consisting
of $n\times n$ matrices such that
\smallskip
\itemitem{(a)}  the entries are either $0$ or $r$th roots of unity,
\smallskip
\itemitem{(b)}  there is exactly one nonzero entry in each row and each 
column.

\smallskip
In the classification [ST] of finite groups generated by complex 
reflections
there is a single infinite family of groups $G(r,p,n)$
and exactly 34 others, the ``exceptional'' complex reflection
groups.  The groups $G(r,p,n)$ are the groups of
$n\times n$ matrices such that
\smallskip
\itemitem{(a)}  the entries are either $0$ or $r$th roots of unity,
\smallskip
\itemitem{(b)}  there is exactly one nonzero entry in each row and each 
column,
\smallskip
\itemitem{(c)}  the $(r/p)$th power of the product of the nonzero entries
is $1$.
\smallskip\noindent
Though we do not know of an early reference which generalizes the theory
of Young tableaux to these groups, it is not difficult to
see that the method that Young uses for the Weyl groups $WD_n$
extends easily to handle the groups $G(r,p,n)$.

Special cases of the groups $G(r,p,n)$ are
\smallskip
\itemitem{(a)}  $G(1,1,n) = S_n$, the symmetric group,
\smallskip
\itemitem{(b)}  $G(2,1,n) = WB_n$, the hyperoctahedral group (i.e.
the Weyl group of type $B_n$),
\smallskip
\itemitem{(c)}  $G(2,2,n) = WD_n$, the Weyl group of type $D_n$,
\smallskip
\itemitem{(d)}  $G(r,1,n) \cong (\ZZ/r\ZZ)\wr S_n= (\ZZ/r\ZZ)^n \rtimes 
S_n$.
\smallskip\noindent
The order of $G(r,1,n)$ is $r^nn!$.
Let $E_{ij}$ be the $n\times n$ matrix with $1$ in the $(i,j)$
position and all other entries $0$.
Then $G(r,1,n)$ can be presented by generators
$$s_1=\zeta E_{11}+\sum_{j\ne 1} E_{ii},
\quad\hbox{and}\quad
s_i = E_{i,i+1} + E_{i+1,i} +\sum_{j\ne i,i+1} E_{jj},
\qquad 2\le i\le n,$$
where $\zeta$ is a primitive $r$th root of unity, and relations
\smallskip
\itemitem{(B1)} $s_is_j = s_js_i$, \hskip.9in if $|i-j|>1$,
\smallskip
\itemitem{(B2)} $s_is_{i+1}s_i = s_{i+1}s_is_{i+1}$,
\qquad for $2\le i\le n-1$,
\smallskip
\itemitem{(BB)} $s_1s_2s_1s_2=s_2s_1s_2s_1$,
\smallskip
\itemitem{(C)} $s_1^r=1$,
\smallskip
\itemitem{(R)} $s_i^2= 1$,\hskip1.2in for $2\le i\le n$.
\smallskip\noindent
The group
$G(r,p,n)$ is the subgroup of index $p$ in $G(r,1,n)$ generated by
$$a_0=s_1^p,\qquad a_1=s_1s_2s_1, \qquad a_i=s_i, \quad 2\le i\le n.$$

\subsection Cyclotomic Hecke algebras  $H_{r,1,n}$.

More recently there has been an interest in Iwahori-Hecke algebras
associated to reflection groups and there has been significant work
generalizing the constructions of A. Young to these algebras.
Iwahori-Hecke algebras of types A, B and D were handled by Hoefsmit [Ho] 
and
other aspects of the theory for these algebras were developed by Dipper,
James and
Murphy [DJ, DJM], Gyoja [Gy] and Wenzl [Wz].  In 1994, Ariki and Koike [AK]
introduced cyclotomic Hecke algebras $H_{r,1,n}$ for the complex reflection
groups $G(r,1,n)$ and they generalized the Young tableau theory to
these algebras.  Theorem 3.18 below shows that the theory of [AK] is
a special case of an even more general theory for affine Hecke algebras.

Let $u_1,\ldots, u_r, q\in \CC$, $q\ne 0$.  The cyclotomic
Hecke algebra $H_{r,1,n}(u_1,\ldots,u_r;q)$ is the algebra over $\CC$ 
given by
generators $T_1,\ldots, T_n$ and relations
\smallskip
\itemitem{(B1)} $T_iT_j = T_jT_i$, \hskip1.8in if $|i-j|>1$,
\smallskip
\itemitem{(B2)} $T_iT_{i+1}T_i = T_{i+1}T_iT_{i+1}$,
\hskip1.15in for $2\le i\le n-1$,
\smallskip
\itemitem{(BB)} $T_1T_2T_1T_2=T_2T_1T_2T_1$,
\smallskip
\itemitem{(qC)} $(T_1-u_1)(T_1-u_2)\cdots (T_1-u_r)=0$,
\smallskip
\itemitem{(qR)} $(T_i-q)(T_i+q^{-1}) = 0$,\hskip1.2in for $2\le i\le n$.
\smallskip\noindent
The algebra $H_{r,1,n}(u_1,\ldots,u_r;q)$
is of dimension $r^n n!$ (see [AK]).
\smallskip
\itemitem{(a)}  $H_{1,1,n}(1;q)$ is the Iwahori-Hecke algebra
of type $A_{n-1}$.
\smallskip
\itemitem{(b)}  $H_{2,1,n}(q,-q^{-1};q)$ is the Iwahori-Hecke algebra
of type $B_n$.
\smallskip
\itemitem{(c)}  If $\zeta$ is a primitive $r$th root of
$1$ then $H_{r,1,n}(1,\zeta,\ldots, \zeta^{r-1};1)$ is the group algebra
$\CC G(r,1,n)$.
\smallskip\noindent
Fact (c) says that the representation theory of the groups $G(r,1,n)$
is a special case of the representation
theory of the algebras $H_{r,1,n}$.

\subsection Cyclotomic Hecke algebras $H_{r,p,n}$.

The work of Brou\'e, Malle and Michel [BMM] demonstrated that there are
cyclotomic Hecke algebras associated to most complex reflection groups
(even exceptional complex reflection groups).   In particular, there
are cyclotomic Hecke algebras $H_{r,p,n}$ corresponding to all the
groups $G(r,p,n)$ and Ariki [Ar2] has generalized the Young tableau
mechanism to these groups (see also [HR]).  Theorems 3.15 and
the mechanism of Theorem 2.8 show that the theory of Ariki is a special 
case of
a general construction for affine Hecke algebras.

Let $r,p,n\in \ZZ_{>0}$
be such that $p$ divides $r$ and let $d=r/p$.
Let $x_0,\ldots, x_{d-1}\in \CC$ and let $\xi$ be a primitive
$p$th root of $1$.
For $1\le j\le r$, define
$$u_j = \xi^kx_\ell,
\qquad \hbox{if $j-1=\ell p+k$,}
\qquad\quad (0\leq k\leq p-1,\  0\leq \ell\leq d-1)$$
i.e., $u_1,\ldots, u_r\in \CC$ are chosen so that
$$(T_1-u_1)(T_1-u_2)\cdots(T_1-u_r)=
(T_1^p-x_0^p)(T_1^p-x_1^p)\cdots(T_1^p-x_{d-1}^p).
$$
The algebra $H_{r,p,n}(x_0,\ldots, x_{d-1};q)$ is the subalgebra of
$H_{r,1,n}(u_1,\ldots, u_r;q)$ generated by the elements
$$
a_0 = T_1^p,\qquad a_1=T_1^{-1}T_2T_1, \qquad
a_i=T_i,\quad \hbox{for $2\le i\le n$.}
$$
Then
\smallskip
\item{(a)}  $H_{2,2,n}(1;q)$ is the Iwahori-Hecke algebra of
type $D_n$,
\smallskip
\item{(b)}  If $\eta$ is a primitive
$d$th root of unity then
$H_{r,p,n}(1,\eta,\cdots,\eta^{d-1};1)$ is the group algebra $\CC G(r,p,n)
$.

\subsection Affine braid groups of type A.

There are three common ways of
depicting affine braids [Cr], [GL], [Jo]:
\smallskip
\item{(a)}  As braids in a (slightly thickened) cylinder,
\smallskip
\item{(b)}  As braids in a (slightly thickened) annulus,
\smallskip
\item{(c)}  As braids with a flagpole.

\smallskip\noindent
See Figure 1.
The multiplication is by placing one cylinder on
top of another,  placing one annulus inside another,
or placing one flagpole braid on top of another.
These are equivalent formulations:
an annulus can be made into a cylinder by turning up the edges,
and a cylindrical braid can be made into a flagpole braid
by putting a flagpole down the middle of the cylinder and
pushing the pole over to the left so that the strings
begin and end to its right.

The group formed by the affine braids with $n$ strands is the
affine braid group of type A.
The affine braid group $\cB_{\infty,1,n}$ is presented by generators
$T_2,\ldots, T_n$ and $X^{\varepsilon_1}$ (see Figure 2) with relations
\smallskip
\itemitem{(B1)} $T_iT_j = T_jT_i$, \hskip1.5in if $|i-j|>1$,
\smallskip
\itemitem{(B2)} $T_iT_{i+1}T_i = T_{i+1}T_iT_{i+1}$,
\hskip.85in for $2\le i\le n-1$,
\smallskip
\itemitem{(BB)} $X^{\varepsilon_1}T_2X^{\varepsilon_1}T_2
=T_2X^{\varepsilon_1}T_2X^{\varepsilon_1}$,
\smallskip
\itemitem{(${\rm B1}'$)} $X^{\varepsilon_1}T_i = T_iX^{\varepsilon_1}$,
\hskip1.3in for $3\le i\le n$.
\smallskip\noindent
Inductively define $X^{\varepsilon_i}\in \cB_{\infty,1,n}$ by
$$X^{\varepsilon_i} = T_i X^{\varepsilon_{i-1}} T_i,
\qquad 2\le i\le n.\formula$$
By drawing pictures of the corresponding affine braids
one can check that the $X^{\varepsilon_i}$ all commute with each
other.  View the symbols $\varepsilon_i$ as a basis of $\RR^n$ so that
$$\RR^n = \sum_{i=1}^n \RR\varepsilon_i,
\qquad\hbox{and let}\qquad
L= \sum_{i=1}^n\ZZ\varepsilon_i.
\formula$$
The affine braid group $\cB_{\infty,1,n}$ contains a large abelian subgroup
$$X = \{ X^\lambda \ |\ \lambda\in L\},
\formula$$
where $X^\lambda=(X^{\varepsilon_1})^{\lambda_1}\cdots
(X^{\varepsilon_n})^{\lambda_n}$
for $\lambda = \lambda_1\varepsilon_1+\cdots+\lambda_n\varepsilon_n\in L$.

The {\it pole winding number} of an affine braid $b\in \cB_{\infty,1,n}$ is
$\kappa(b)$ where
$\kappa\colon \cB_{\infty,1,n} \to \ZZ$
is the group homomorphism defined by
$\kappa(X^{\varepsilon_1})=1$ and $\kappa(T_i)=0,\quad 2\le i\le n.$
The affine braid group $\cB_{\infty,p,n}$
is the subgroup of $\cB_{\infty,1,n}$ of affine
braids with pole winding number equal to $0$ (mod $p$),
$$\cB_{\infty,p,n} =
\{ b\in \cB_{\infty,1,n} \ |\ \kappa(b)=0 {\rm\ \  (mod\ } p)\}.
\formula$$
Define
$$
Q=\sum_{i=2}^n \ZZ(\varepsilon_i-\varepsilon_{i-1})
\qquad\hbox{and}\qquad
L_p = Q+\sum_{i=1}^n p\ZZ\varepsilon_i,
\formula$$
for each nonnegative integer $p$.   The lattice $L_p$ is a lattice of index
$p$ in $L$.
Then
$$\cB_{\infty,p,n} = \langle X^\lambda, T_i\
|\ \lambda\in L_p,  2\le i\le n\rangle$$
and the group $X^{L_p}=\langle X^\lambda\ |\ \lambda\in L_p\rangle$ is
an abelian subgroup of $\cB_{\infty,p,n}$.

\subsection Affine Hecke algebras of type A.

The affine Hecke algebra $H_{\infty,1,n}$
(resp. $H_{\infty,p,n}$)
is the quotient of the group algebra $\CC\cB_{\infty,1,n}$
(resp. $\CC\cB_{\infty,p,n}$)
by the relations
$$T_i^2 = (q-q^{-1})T_i+1, \qquad 2\le i\le n.$$
Let $L$ and $X^\lambda$ be as in (1.5) and (1.6).  The subalgebra
$$\CC[X] = {\rm span}\{ X^\lambda\ |\ \lambda\in L\}\qquad
(\hbox{resp.}\ \   \CC[X^{L_p}] = {\rm span}\{ X^\lambda\ |\ \lambda\in 
L_p\})
\formula$$
is a commutative subalgebra of $H_{\infty,1,n}$ (resp. $H_{\infty,p,n}$).
The symmetric group $S_n$ acts on the lattice $L$ by permuting
the $\varepsilon_i$ and the lattices $L_p$ are $S_n$-invariant sublattices 
of
$L$.  Let $s_i=(i,i-1)\in S_n$ and 
$\alpha_i=\varepsilon_i-\varepsilon_{i-1}$.
For $2\le i\le n$ and $\lambda\in L$ (resp. $\lambda\in L_p$),
$$
X^\lambda T_i
= X^{s_i\lambda}T_i + (q-q^{-1})
{X^\lambda -X^{s_i\lambda}\over 1-X^{-\alpha_i}},
\formula
$$
as elements of $H_{\infty,1,n}$ (resp. $H_{\infty,p,n}$).

For each element $w\in S_n$ define $T_w = T_{i_1}\cdots T_{i_p}$
if $w=s_{i_1}\cdots s_{i_p}$ is an expression of $w$ as a product of
simple reflections $s_i$ such that $p$ is minimal.  The element $T_w$ does 
not
depend on the choice of the reduced expression of $w$ [Bou, Ch. IV \S 2 Ex.
23].
The sets
$$\{X^\lambda T_w \ |\ \lambda\in L, w\in S_n\}
\qquad\hbox{and}\qquad
\{X^\lambda T_w \ |\ \lambda\in L_p, w\in S_n\}$$
are bases of $H_{\infty,1,n}$ and $H_{\infty,p,n}$, respectively [Lu].
The center of $H_{\infty,1,n}$ is
$$Z(H_{\infty,1,n}) = \CC[X^{\varepsilon_1},\ldots,X^{\varepsilon_n}]^{S_n}
=\CC[X]^{S_n},
\formula$$
and $\CC[X^{L_p}]^{S_n}$ is the center of $H_{\infty,p,n}$ (see Theorem
4.12 below).

\section 2.  Representation theory transfer

In this section we provide the mechanism for obtaining the representation
theory of $H_{\infty,p,n}$ from $H_{\infty,1,n}$ and for obtaining the
representation theory of cyclotomic Hecke algebras from affine Hecke 
algebras.
In order to obtain the representation theory of $H_{\infty,p,n}$ from
$H_{\infty,1,n}$ we identify $H_{\infty,p,n}$ as the set of fixed points
of a certain group $G$ acting on $H_{\infty,1,n}$ by automorphisms.  Once
this is done, the extended version of Clifford theory given in the Appendix
allows one to construct the representations of $H_{\infty,p,n}$ from
those of $H_{\infty,1,n}$.  The same technique can be applied to obtain the
representations of the braid groups ${\cal B}_{\infty,p,n}$ from those
of ${\cal B}_{\infty,1,n}$, of the complex reflection groups $G(r,p,n)$
from those of $G(r,1,n)$, and of the Weyl groups $WD_n$ from those of the
Weyl groups $WB_n$.

\subsection  Obtaining $H_{\infty,p,n}$-modules from 
$H_{\infty,1,n}$-modules.

The following result is what is needed to apply the Clifford theory
developed in
the Appendix to derive the representation theory of the algebras
$H_{\infty,p,n}$ from that of $H_{\infty,1,n}$.

\thm
Let $\xi$ be a primitive $p$th root of unity.  The algebra automorphism
$g\colon H_{\infty,1,n}\to H_{\infty,1,n}$ defined by
$$g(X^{\varepsilon_1}) = \xi X^{\varepsilon_1},
\qquad\hbox{and}\qquad g(T_i)=T_i, \quad 2\le i\le n,$$
gives rise to an action of the group
$\ZZ/p\ZZ = \{1,g,\ldots, g^{p-1}\}$ on $H_{\infty,1,n}$ by
algebra automorphisms and
$$H_{\infty,p,n} = (H_{\infty,1,n})^{\ZZ/p\ZZ}$$
is the set of fixed points of the $\ZZ/p\ZZ$-action.
\pf
Immediate from the definitions of $H_{\infty,1,n}$, $H_{\infty,p,n}$
and (1.6).
\endpf

The action of $\ZZ/p\ZZ$ on $H_{\infty,1,n}$ which is given in
Theorem 2.2 induces an action of $\ZZ/p\ZZ$ on the simple
$H_{\infty,1,n}$-modules, see (A.1) in the Appendix.  The stabilizer $K$ 
of the
action of $\ZZ/p\ZZ$ on a simple $H_{\infty,1,n}$-module $M$ is the
{\it inertia group} $K$ of $M$.  The action of
$K$ commutes with the action of $H_{\infty,p,n}$ on $M$
and we have a decomposition
$$M\cong \bigoplus_{j=0}^{|K|-1} M^{(j)}\otimes K^{(j)},
\formula$$
where $K^{(j)}$, $1\le j\le |K|-1$, are the simple $K$-modules
and $M^{(j)}$ are $H_{\infty,p,n}$-modules.  Theorem
A.13 of the Appendix shows that the $M^{(j)}$ are simple
$H_{\infty,p,n}$-modules and that all simple $H_{\infty,p,n}$-modules
are constructed in this way.  In Theorem 3.15
we  show that this method gives a combinatorial construction of the
module
$M^{(j)}$ in any case when the Young tableau theory is available.
This is a generalization of the method used in [Ar2] and [HR].

\subsection The surjective algebra homomorphisms $\Phi$ and $\Phi_p$.

The homomorphisms $\Phi$ and $\Phi_p$
described below are the primary tools for
transferring results from the affine Hecke algebras to
cyclotomic Hecke algebras.  Many results are easier to prove
for affine Hecke algebras because of the large commutative
subalgebra $\CC[X]$ which is available in the affine Hecke
algebra.  The homomorphism $\Phi$ has also been used
by Cherednik [Ch2], Ariki [Ar3] and many others.

\prop  Let $H_{\infty,1,n}$ be the affine Hecke algebra of type A defined 
in
(1.9) and let $H_{r,1,n}(u_1,\ldots,u_r;q)$ denote the cyclotomic Hecke 
algebra
of (1.1).
\smallskip\noindent
(a) Fix $u_1,\ldots, u_r\in \CC$, $q\in \CC^*$.
There is a surjective algebra homomorphism given by
$$\matrix{
\Phi\colon &H_{\infty,1,n} &\longrightarrow &H_{r,1,n}(u_1,\ldots, u_r;q) 
\cr
&T_i &\longmapsto &T_i, &\qquad &2\le i\le n,\cr
&X^{\varepsilon_1} &\longmapsto &T_1. \cr}$$
(b) Restricting the homomorphism $\Phi$ to $H_{\infty,p,n}$ yields a 
surjective
homomorphism defined by
$$\matrix{
\Phi_p\colon &H_{\infty,p,n} &\longrightarrow
&H_{r,p,n}(x_0,\ldots,x_{d-1};q) \cr
&T_i &\longmapsto &a_i, &\qquad &2\le i\le n,\cr
&X^{p\varepsilon_1} &\longmapsto &a_0, \cr
&X^{\varepsilon_2-\varepsilon_1} &\longmapsto &a_1a_2. \cr}$$
\pf
The result
follows directly from the definitions of the affine Hecke
algebras $H_{\infty,1,n}$ and $H_{\infty,p,n}$
(see (1.9)) and the cyclotomic
Hecke algebras $H_{r,1,n}$ and $H_{r,p,n}$ (see (1.1) and (1.2)).
\endpf

Let $L_p$ be the lattice defined in (1.8) and define
$C_r=\{ \lambda=\sum_{i=1}^n
\lambda_i\varepsilon_i \ |\  0\le \lambda_i< r \}.$
Ariki and Koike [AK] and Ariki [Ar2] have shown that the sets
$$\{\Phi(X^\lambda T_w) \ |\ \lambda\in C_r, w\in S_n\}
\qquad\hbox{and}\qquad
\{\Phi_p(X^\lambda T_w) \ |\ \lambda\in C_r\cap L_p, w\in S_n\}
\formula$$
are bases of $H_{r,1,n}$ and $H_{r,p,n}$, respectively.

\subsection Relating $H_{\infty,1,n}$-modules and $H_{r,1,n}$-modules.

The representation theory of
the affine Hecke algebra $H_{\infty,1,n}$ is equivalent to the
representation theory of the cyclotomic Hecke algebras $H_{r,1,n}$
(considering all possible $u_1,\ldots, u_r\in \CC$).
The elementary constructions in the following theorem allow us to make
$H_{\infty,1,n}$-modules into $H_{r,1,n}$-modules {\it and vice versa}.

\thm  Let $H_{\infty,1,n}$ be the affine Hecke algebra of type A defined in
(1.9) and let $H_{r,1,n}(u_1,\ldots,u_r;q)$ be the cyclotomic Hecke algebra
of (1.1).
\item{(a)}
Fix $u_1,\ldots, u_r\in \CC$, $q\in \CC^*$, and let $\Phi$
be the surjective homomorphism of Proposition 2.5.
If $M$ is a simple $H_{r,1,n}(u_1,\ldots, u_r;q)$-module then
defining
$$hm = \Phi(h)m,
\qquad\hbox{for all $h\in H_{\infty,1,n}$ and all $m\in M$,}$$
makes $M$ into a simple $H_{\infty,1,n}$-module.
\smallskip
\item{(b)} Let $M$ be a simple $H_{\infty,1,n}$-module and let
$\rho\colon H_{\infty,1,n}\to \End(M)$ be the corresponding
representation.  Let $u_1,\ldots, u_r\in \CC$ be such that
the minimal polynomial $p(t)$ of the matrix
$\rho(X^{\varepsilon_1})$
divides the polynomial $(t-u_1)\cdots (t-u_r)$.
Define an action of $H_{r,1,n}(u_1,\ldots,u_r;q)$ on $M$ by
$$
T_1m = X^{\varepsilon_1}m
\qquad\hbox{and}\qquad T_im = T_im, \quad\hbox{$2\le i\le n$,}
$$
for all $m\in M$.
Then $M$ is a simple $H_{r,1,n}(u_1,\ldots,u_r;q)$-module.
\pf
The Theorem follows directly from the definitions of
$H_{\infty,1,n}$ and $H_{r,1,n}(u_1,\ldots, u_r;q)$, and
the construction of the surjective homomorphism $\Phi$.
\endpf

\remark
The same translations work for arbitrary finite dimensional
modules; in particular, they work for indecomposable modules
and {\it preserve composition series}.
\endprop

\section 3.  Standard Young tableaux, representations and Jucys-Murphy 
elements

In this section we review the generalization of standard Young tableaux
in [Ra5] which is used to construct representations of the affine Hecke
algebras $H_{\infty,1,n}$. Then we show how this theory can be transported 
to
provide combinatorial constructions of simple modules for the cyclotomic
Hecke algebras $H_{r,1,n}(u_1,\ldots,u_r;q)$ and
$H_{r,p,n}(x_0,\ldots,x_{d-1};q)$.
This approach shows that Jucy-Murphy type elements in the
cyclotomic Hecke algebras
arise naturally as images of the elements $X^{\varepsilon_i}$ in
the affine Hecke algebra.  The standard Jucys-Murphy type theorems then 
follow
almost immediately from standard affine Hecke algebra facts.

\subsection{Skew shapes and standard tableaux.}

A partition $\lambda$ is a collection of $n$
boxes in a  corner.  We shall conform to the conventions in [Mac] and 
assume
that gravity goes up and to the left.
$$
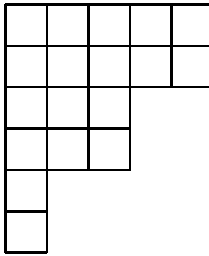
\beginpicture
\setcoordinatesystem units <0.5cm,0.5cm>         % sets scale
\setplotarea x from 0 to 4, y from 0 to 3    % sets plot size up
\linethickness=0.5pt                          % sets line thickness
\putrule from 0 6 to 5 6          %
\putrule from 0 5 to 5 5          %  draws horizontal lines
\putrule from 0 4 to 5 4          %
\putrule from 0 3 to 3 3          %
\putrule from 0 2 to 3 2          %
\putrule from 0 1 to 1 1          %
\putrule from 0 0 to 1 0          %

\putrule from 0 0 to 0 6        %
\putrule from 1 0 to 1 6        %
\putrule from 2 2 to 2 6        %
\putrule from 3 2 to 3 6        %  draws vertical lines
\putrule from 4 4 to 4 6        %
\putrule from 5 4 to 5 6        %
%\vshade 2 1 2   3 1 2 /
\endpicture
$$
Any partition $\lambda$ can be identified with
the sequence $\lambda=(\lambda_1\ge \lambda_2\ge \ldots )$
where $\lambda_i$ is the number of boxes in row $i$ of $\lambda$.
The rows and columns are numbered in the same way as for matrices.
In the example above we have $\lambda=(553311)$.
If $\lambda$ and $\mu$ are partitions such that $\mu_i\le \lambda_i$
for all $i$ we write $\mu\subseteq \lambda$.
The {\it skew shape} $\lambda/\mu$
consists of all boxes of $\lambda$ which are not in $\mu$.
Any skew shape is a union of connected components.  Number the boxes of 
each
skew shape $\lambda/\mu$ along major diagonals from southwest to northeast 
and
$$\hbox{write ${\rm box}_i$ to indicate the box numbered $i$.}$$

Let $\lambda/\mu$ be a skew shape with $n$ boxes.
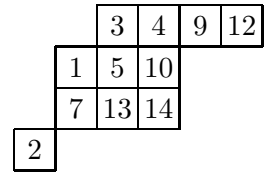
A {\it standard tableau of shape} $\lambda/\mu$ is a filling
of the boxes in the skew shape
$\lambda/\mu$ with the numbers $1,\ldots,n$ such that
the numbers increase from left to right in each row and from top to
bottom down each column.

\subsection{Placed skew shapes.}

Let $\RR+i[0,2\pi/\ln(q^2))=\{ a+bi\ |\ a\in \RR, 0\le b\le 2\pi/\ln(q^2)
\}\subseteq \CC$.  If $q$ is a positive real number then
the function
$$\matrix{
\RR+i[0,2\pi/\ln(q^2)) &\longrightarrow &\CC^*\hfill \cr
x &\longmapsto &q^{2x}=e^{\ln(q^2)x} \cr}$$
is a bijection.
The elements of $[0,1)+i[0,2\pi/\ln(q^2))$ index the $\ZZ$-cosets in
$\RR+i[0,2\pi/\ln(q^2))$.

A {\it placed skew shape} is a pair $(c,\lambda/\mu)$
consisting of a skew shape $\lambda/\mu$ and a {\it content
function}
$$c\colon \{\hbox{boxes of $\lambda/\mu$}\} \longrightarrow
\RR+i[0,2\pi/\ln(q^2))
\qquad\hbox{such that}\formula$$
$$
\matrix{
c({\rm box}_j)- c({\rm box}_i)\ge 0, \hfill
&\quad &\hbox{if $i<j$ and $c({\rm box}_j)-c({\rm box}_i)\in \ZZ$}, \cr
c({\rm box}_j)=c({\rm box}_i)+1,\hfill
&&\hbox{if and only if ${\rm box}_i$ and ${\rm box}_j$
are on adjacent diagonals, and} \cr
c({\rm box}_i)=c({\rm box}_j),\hfill
&&\hbox{if and only if ${\rm box}_i$ and ${\rm box}_j$
are on the same diagonal.} \cr
}$$
This is a generalization of the usual notion of the content of a box
in a partition (see [Mac] I \S 1 Ex. 3).

Suppose that $(c,\lambda/\mu)$ is a placed skew shape such that
$c$ takes values in $\ZZ$.  One can visualize $(c,\lambda/\mu)$
by placing
$\lambda/\mu$ on a piece of infinite graph paper
where the diagonals of the
graph paper are indexed consecutively (with elements of $\ZZ$) from 
southeast
to northwest.  The {\it content} of a box $b$ is the index $c(b)$ of the
diagonal that $b$ is on.
In the general case, when $c$ takes values in
$\RR+i[0,2\pi/\ln(q^2))$, one imagines
a book where each page is a sheet of infinite graph paper
with the diagonals indexed consecutively
(with elements of $\ZZ$) from southeast to northwest.
The pages are numbered by values $\beta\in [0,1)+i[0,2\pi/\ln(q^2))$
and there is a skew shape $\lambda^{(\beta)}/\mu^{(\beta)}$ placed on page
$\beta$.    The skew shape $\lambda/\mu$ is a union of the disjoint skew 
shapes
on each page,
$$\lambda/\mu = \bigsqcup_\beta 
\left(\lambda^{(\beta)}/\mu^{(\beta)}\right),
\qquad \beta\in [0,1)+i[0,2\pi/\ln(q^2)),\formula$$
and the content function is given by
$$
c(b) =
\hbox{(page number of the page containing $b$)}
+ \hbox{(index of the diagonal containing $b$).}
\formula$$
for a box $b\in \lambda/\mu$.

\subsection{Example.}

The following diagrams illustrate standard tableaux and
the numbering of boxes in a skew shape $\lambda/\mu$.
$$
\matrix{
\beginpicture
\setcoordinatesystem units <0.5cm,0.5cm>         % sets scale
\setplotarea x from 0 to 4, y from 0 to 3    % sets plot size up
\linethickness=0.5pt                          % sets line thickness
\put{1} at 0.5 0.5
\put{2} at 0.5 1.5
\put{3} at 1.5 1.5
\put{4}  at 3.5 2.5
\put{5} at 4.5 3.5
\put{6} at 4.5 4.5
\put{7}  at 5.5 3.5
\put{8}  at 5.5 4.5
\put{10}  at 5.5 5.5
\put{9}  at 6.5 3.5
\put{11}  at 6.5 4.5
\put{12}  at 6.5 5.5
\put{13}  at 7.5 5.5
\put{14}  at 8.5 5.5
\putrule from 5 6 to 9 6          %
\putrule from 4 5 to 9 5          %  draws horizontal lines
\putrule from 4 4 to 7 4          %
\putrule from 3 3 to 7 3          %
\putrule from 3 2 to 4 2          %
\putrule from 0 2 to 2 2          %
\putrule from 0 1 to 2 1          %
\putrule from 0 0 to 1 0          %
\putrule from 0 0 to 0 2        %
\putrule from 1 0 to 1 2        %
\putrule from 2 1 to 2 2        %
\putrule from 3 2 to 3 3        %  draws vertical lines
\putrule from 4 2 to 4 5        %
\putrule from 5 3 to 5 6        %
\putrule from 6 3 to 6 6        %
\putrule from 7 3 to 7 6        %
\putrule from 8 5 to 8 6        %
\putrule from 9 5 to 9 6        %
%\vshade 2 1 2   3 1 2 /
\endpicture
&\qquad\qquad
&
\beginpicture
\setcoordinatesystem units <0.5cm,0.5cm>         % sets scale
\setplotarea x from 0 to 4, y from 0 to 3    % sets plot size up
\linethickness=0.5pt                          % sets line thickness
\put{11} at 0.5 0.5
\put{6} at 0.5 1.5
\put{8} at 1.5 1.5
\put{2}  at 3.5 2.5
\put{7} at 4.5 3.5
\put{1} at 4.5 4.5
\put{13}  at 5.5 3.5
\put{5}  at 5.5 4.5
\put{3}  at 5.5 5.5
\put{14}  at 6.5 3.5
\put{10}  at 6.5 4.5
\put{4}  at 6.5 5.5
\put{9}  at 7.5 5.5
\put{12}  at 8.5 5.5
\putrule from 5 6 to 9 6          %
\putrule from 4 5 to 9 5          %  draws horizontal lines
\putrule from 4 4 to 7 4          %
\putrule from 3 3 to 7 3          %
\putrule from 3 2 to 4 2          %
\putrule from 0 2 to 2 2          %
\putrule from 0 1 to 2 1          %
\putrule from 0 0 to 1 0          %
\putrule from 0 0 to 0 2        %
\putrule from 1 0 to 1 2        %
\putrule from 2 1 to 2 2        %
\putrule from 3 2 to 3 3        %  draws vertical lines
\putrule from 4 2 to 4 5        %
\putrule from 5 3 to 5 6        %
\putrule from 6 3 to 6 6        %
\putrule from 7 3 to 7 6        %
\putrule from 8 5 to 8 6        %
\putrule from 9 5 to 9 6        %
%\vshade 2 1 2   3 1 2 /
\endpicture
\cr
\hbox{$\lambda/\mu$ with boxes numbered}
&&\hbox{A standard tableau $L$ of shape $\lambda/\mu$} \cr
}
$$

The following picture shows the contents of the boxes in the
placed skew shape $(c,\lambda/\mu)$ such that the sequence
$(c({\rm box}_1),\ldots, c({\rm box}_n))$ is
$(-7,-6,-5,-2,0,1,1,2,2,3,3,4,5,6)$.
$$
\matrix{
\beginpicture
\setcoordinatesystem units <0.5cm,0.5cm>         % sets scale
\setplotarea x from 0 to 4, y from 0 to 3    % sets plot size up
\linethickness=0.5pt                          % sets line thickness
\put{-7} at 0.5 0.5
\put{-6} at 0.5 1.5
\put{-5} at 1.5 1.5
\put{-2}  at 3.5 2.5
\put{0} at 4.5 3.5
\put{1} at 4.5 4.5
\put{1}  at 5.5 3.5
\put{2}  at 5.5 4.5
\put{3}  at 5.5 5.5
\put{2}  at 6.5 3.5
\put{3}  at 6.5 4.5
\put{4}  at 6.5 5.5
\put{5}  at 7.5 5.5
\put{6}  at 8.5 5.5
\putrule from 5 6 to 9 6          %
\putrule from 4 5 to 9 5          %  draws horizontal lines
\putrule from 4 4 to 7 4          %
\putrule from 3 3 to 7 3          %
\putrule from 3 2 to 4 2          %
\putrule from 0 2 to 2 2          %
\putrule from 0 1 to 2 1          %
\putrule from 0 0 to 1 0          %
\putrule from 0 0 to 0 2        %
\putrule from 1 0 to 1 2        %
\putrule from 2 1 to 2 2        %
\putrule from 3 2 to 3 3        %  draws vertical lines
\putrule from 4 2 to 4 5        %
\putrule from 5 3 to 5 6        %
\putrule from 6 3 to 6 6        %
\putrule from 7 3 to 7 6        %
\putrule from 8 5 to 8 6        %
\putrule from 9 5 to 9 6        %
%\vshade 2 1 2   3 1 2 /
\endpicture
\cr
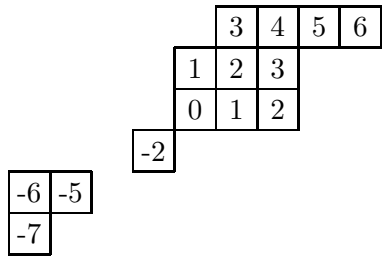
\hbox{Contents of the boxes of $(c,\lambda/\mu)$} \cr
}
$$
The following picture shows the contents of the boxes in the placed skew
shape $(c,\lambda/\mu)$ such that
$(c({\rm box}_1),\ldots, c({\rm box}_n))
=(-7,-6,-5,-3/2,1/2,3/2,3/2,5/2,5/2,7/2,7/2,9/2,11/2,13/2)$.
$$
\beginpicture
\setcoordinatesystem units <0.5cm,0.5cm>         % sets scale
\setplotarea x from -4 to 4, y from 0 to 7    % sets plot size up
\linethickness=0.5pt                          % sets line thickness
\put{-7} at -3.5 3.5
\put{-6} at -3.5 4.5
\put{-5} at -2.5 4.5
\put{-${3\over2}$}  at 3.5 2.5
\put{$1\over2$} at 4.5 3.5
\put{$3\over2$} at 4.5 4.5
\put{$3\over2$}  at 5.5 3.5
\put{$5\over 2$}  at 5.5 4.5
\put{$7\over 2$}  at 5.5 5.5
\put{$5\over 2$}  at 6.5 3.5
\put{$7\over2$}  at 6.5 4.5
\put{$9\over2$}  at 6.5 5.5
\put{$11\over2$}  at 7.5 5.5
\put{$13\over2$}  at 8.5 5.5
\putrule from -4 5 to -2 5          %
\putrule from -4 4 to -2 4          %   horizontal lines on first page
\putrule from -4 3 to -3 3          %
\putrule from -4 3 to -4 5        %
\putrule from -3 3 to -3 5        %     Vertical lines on first page
\putrule from -2 4 to -2 5        %
\putrule from 5 6 to 9 6          %
\putrule from 4 5 to 9 5          %  draws horizontal lines
\putrule from 4 4 to 7 4          %   on second page
\putrule from 3 3 to 7 3          %
\putrule from 3 2 to 4 2          %
\putrule from 3 2 to 3 3        %  draws vertical lines
\putrule from 4 2 to 4 5        %    on second page
\putrule from 5 3 to 5 6        %
\putrule from 6 3 to 6 6        %
\putrule from 7 3 to 7 6        %
\putrule from 8 5 to 8 6        %
\putrule from 9 5 to 9 6        %
\setdashes
\putrule from 0 0 to 0 7           % draws break between pages
\put{0}   at -3.5 0.5      %
\put{1/2} at 5 0.5       % labels pages
\endpicture
$$
This ``book'' has two pages, with page numbers $0$ and $1/2$.
\endpf

\subsection  Calibrated $H_{\infty,1,n}$-modules.

A finite dimensional $\tilde H_{\infty,1,n}$-module
$M$ is {\it calibrated} if it has a basis
$\{v_t\}$ such that for each $\lambda\in L$ and each $v_t$ in the basis
$$X^\lambda v_t = t(X^\lambda)v_t,\qquad
\hbox{for some $t(X^\lambda)\in \CC$.}$$
This is the class of representations of the affine Hecke algebra for which
there
is a good theory of Young tableaux [Ra2].

The following theorem classifies and constructs all irreducible calibrated
representations of the affine Hecke algebra $H_{\infty,1,n}$.  The 
construction
is a direct generalization of A. Young's classical
``seminormal construction'' of the irreducible representations of the 
symmetric
group [Yg2]. Young's construction was generalized to Iwahori-Hecke 
algebras of
type A by Hoefsmit [Hf] and Wenzl [Wz] independently, to Iwahori-Hecke 
algebras
of types B and D by Hoefsmit [Hf] and to cyclotomic Hecke algebras by 
Ariki and
Koike [AK].   In (3.7) and (3.11) below we
  show how all of these earlier
generalizations of Young's construction can be obtained from Theorem 3.8.  
Some
parts of Theorem 3.8 are originally due to I. Cherednik, and are stated
in [Ch1, \S 3].

\thm {\rm ([Ra, Theorem 4.1])} Let $(c,\lambda/\mu)$ be a placed skew shape
with
$n$ boxes.  Define an action of $H_{\infty,1,n}$ on the vector space
$$\tilde H^{(c,\lambda/\mu)} =
\CC{\rm -span}\{ v_L\ |\ \hbox{$L$ is a standard
tableau of shape $\lambda/\mu$}\}$$
by the formulas
$$\eqalign{
X^{\varepsilon_i} v_L &= q^{2c(L(i))} v_L, \cr
T_i v_L &= (T_i)_{LL} v_L + (q^{-1}+(T_i)_{LL}) v_{s_iL}, \cr
}
$$
where $s_iL$ is the same as $L$ except that the entries $i-1$ and $i$ are
interchanged,
$$(T_i)_{LL} = { q-q^{-1} \over 1-q^{2(c(L(i-1))-c(L(i)))} },
\quad\quad
\hbox{$v_{s_iL}=0$ if $s_iL$ is not a standard tableau,}$$
and $L(i)$ denotes the box of $L$ containing the entry $i$.
\smallskip
\item{(a)}  $\tilde H^{(c,\lambda/\mu)}$ is a calibrated irreducible
$H_{\infty,1,n}$-module.
\smallskip
\item{(b)}  The modules $\tilde H^{(c,\lambda/\mu)}$ are
non-isomorphic.
\smallskip
\item{(c)}  Every irreducible calibrated $H_{\infty,1,n}$-module is
isomorphic to $\tilde H^{(c,\lambda/\mu)}$ for some placed skew shape
$(c,\lambda/\mu)$.
\endthm

\remark  All of the irreducible modules for the affine Hecke algebra
have been classified and constructed by Kazhdan and Lusztig [KL].
The construction in [KL] is geometric and noncombinatorial.  It is
nontrivial to relate the construction of Theorem 3.8 and
the classification in [KL].   This point is the subject of [RR].
\endprop

\subsection Calibrated $H_{\infty,p,n}$-modules.

A finite dimensional $H_{\infty,p,n}$ module
$M$ is {\it calibrated} if it has a basis
$\{v_t\}$ such that for each $\lambda\in L_p$ and each $v_t$ in the basis
$$X^\lambda v_t = t(X^\lambda)v_t,\qquad
\hbox{for some $t(X^\lambda)\in \CC$.}$$
Let us show how Theorem 3.8, Theorem 2.2 and Theorem A.13 provide explicit
constructions of simple calibrated $H_{\infty,p,n}$-modules.
The resulting construction is a generalization of the construction
of $H_{r,p,n}(x_0,\ldots,x_{d-1};q)$-modules given by Ariki [Ar2]
(as amplified and applied in [HR]).  Comparing the following
machinations with those in [HR, \S 3] (where more pictures are given) will 
be
helpful.

The $(\ZZ/p\ZZ)$-action on $H_{\infty,1,n}$ induces an action
of $\ZZ/p\ZZ$ on the simple $H_{\infty,1,n}$-modules, as in (A.1) of the
Appendix,
and this action takes calibrated modules
to calibrated modules since the $(\ZZ/p\ZZ)$-action on
$H_{\infty,1,n}$ preserves the subalgebra $\CC[X]$.  The $\ZZ/p\ZZ$
action on simple calibrated modules can be described combinatorially as
follows.

If $(c,\lambda/\mu)$ is a placed skew shape with $n$ boxes and
$g\in \ZZ/p\ZZ$ define
$$g(c,\lambda/\mu)=(c-i\alpha/p,\lambda/\mu),
\qquad\hbox{where $\alpha = 2\pi/\ln(q^2)$,}
\formula$$
and $c-i\alpha/p$ denotes the content function defined by
$(c-i\alpha/p)(b) = c(b)-i\alpha/p$, for all boxes $b\in \lambda/\mu$.
To make this definition we are identifying the set $[0,\alpha)$
with $\RR/\alpha\ZZ$. One can imagine the placed skew shape as a book
with pages numbered by values
$\beta\in [0,1)+i(\RR/\alpha\ZZ)$ and a skew shape
$\lambda^{(\beta)}/\mu^{(\beta)}$ on each page.  The action of
$$\hbox{$g$ cyclically permutes the pages numbered
$\beta+i(k/p)\alpha$,\quad $0\le k< p$.}$$
If $L$ is a standard tableau of shape
$(c,\lambda/\mu)$ let
$gL$ denote the same filling of $\lambda/\mu$ as $L$
but viewed as a standard tableaux of shape $g(c,\lambda/\mu)$.

Let $^g\tilde H^{(c,\lambda/\mu)}$ be the $H_{\infty,1,n}$-module
$\tilde H^{(c,\lambda/\mu)}$ except twisted by the automorphism $g$,
see (A.1) in the Appendix.
It follows from the formulas in Theorem 3.8 that the map
$$\matrix{
\phi\colon\  &^g{\tilde H}^{(c,\lambda/\mu)} &\longrightarrow
&\tilde H^{g(c,\lambda/\mu)} \cr
&v_L &\longmapsto &v_{gL}\cr}
\formula$$
is an $H_{\infty,1,n}$-module isomorphism.
Indeed, since $g^{-1}(T_j)=T_j$ and
$(T_j)_{gL,gL}=(T_j)_{LL}$,
$$\eqalign{
\phi(T_j\circ v_L)&=T_j\phi(v_L), \quad\hbox{and} \cr
\phi(X^{\varepsilon_j}\circ v_L)
&=\phi(g^{-1}(X^{\varepsilon_j})v_L)
=\phi(\xi^{-1} X^{\varepsilon_j}v_L)
=q^{2(c(L(j))-i\alpha/p)}v_{gL}
=X^{\varepsilon_j}v_{gL}
=X^{\varepsilon_j}\phi(v_L),\cr}$$
where $\circ$ denotes the $H_{\infty,1,n}$-action on
$^g\tilde H^{(c,\lambda/\mu)}$ as in (A.1).
Identify $^g{\tilde H}^{(c,\lambda/\mu)}$ with
$\tilde H^{g(c,\lambda/\mu)}$ via the isomorphism in (3.12).

Let $(c,\lambda/\mu)$ be a placed skew shape with $n$ boxes and let
$K_{(c,\lambda/\mu)}$ be the stabilizer of $(c,\lambda/\mu)$ under the
action of
$\ZZ/p\ZZ$.  The cyclic group $K_{(c,\lambda/\mu)}$ is a realization of the
inertia group of
$\tilde H^{(c,\lambda/\mu)}$ and
$$K_{(c,\lambda/\mu)}=\{(g^{\kappa})^\ell\colon
\tilde H^{(c,\lambda/\mu)}\to \tilde H^{(c,\lambda/\mu)}
\ |\ 0\le \ell \le |K_{c,\lambda/\mu)}|-1\},$$
where $\kappa$ is the smallest integer between $1$ and $p$
such that $g^\kappa(c,\lambda/\mu)=(c,\lambda/\mu)$ and 
$|K_{(c,\lambda/\mu)}|$
is the order of $K_{(c,\lambda/\mu)}$.  The elements of
$K_{(c,\lambda/\mu)}$ are all $H_{\infty,p,n}$-module isomorphisms.
Since $K_{(c,\lambda/\mu)}$ is a cyclic group the
irreducible $K_{(c,\lambda/\mu)}$-modules are all one-dimensional and the
characters of these modules are given explicitly by
$$\matrix{
\eta_j\colon &K_{(c,\lambda/\mu)} &\longrightarrow & \CC \cr
&g^\kappa &\longmapsto &\xi^{j\kappa}, \cr}
\qquad \qquad 0\le j\le |K_{(c,\lambda/\mu)}|-1,$$
since $\xi^\kappa$ is a primitive $|K_{(c,\lambda/\mu)}|$-th root of unity.
The element
$$p_j = \sum_{\ell=0}^{|K_{(c,\lambda/\mu)}|-1}
\xi^{-j\ell\kappa}g^{\ell\kappa}.
\formula$$
is the minimal idempotent of the group algebra $\CC
K_{(c,\lambda/\mu)}$ corresponding to the irreducible character $\eta_j$.
It follows (from a standard double centralizer result, [Bou2]) that, as an
$(H_{\infty,p,n}, K_{(c,\lambda/\mu)})$-bimodule,
$$\tilde H^{(c,\lambda/\mu)} \cong
\bigoplus_{j=0}^{|K_{(c,\lambda/\mu)}|-1}
\tilde H^{(c,\lambda/\mu,j)}\otimes K^{(j)},
\qquad\hbox{where}\quad
\tilde H^{(c,\lambda/\mu,j)} = p_j\tilde H^{(c,\lambda/\mu)},
\formula$$
and $K^{(j)}$ is the irreducible $K_{(c,\lambda/\mu)}$-module with 
character
$\eta_j$.  The following theorem now follows from Theorem A.13 of the 
Appendix.

\thm
Let $(c,\lambda/\mu)$ be a placed skew shape with $n$ boxes and
let $\tilde
H^{(c,\lambda/\mu)}$ be the simple calibrated $H_{\infty,1,n}$-module
constructed in Theorem 3.8.  Let $K_{(c,\lambda/\mu)}$ be the inertia
group of $\tilde H^{(c,\lambda/\mu)}$ corresponding to the action of
$\ZZ/p\ZZ$ on $H_{\infty,1,n}$ defined by (A.1).  If $p_j$ is the minimal
idempotent of $K_{(c,\lambda/\mu)}$ given by (3.13) then
$$\tilde H^{(c,\lambda/\mu,j)} = p_j\tilde H^{(c,\lambda/\mu)}$$
is a simple calibrated $H_{\infty,p,n}$-module.
\endthm

Theorem 3.15 provides a generalization of the construction of
the $H_{r,p,n}(x_0,\ldots, x_{d-1};q)$-modules which was given by
Ariki [Ar2] and extended and applied in [HR].

\subsection Simple $H_{r,1,n}(u_1,\ldots,u_r;q)$-modules.

Many (usually all) of the simple $H_{r,1,n}(u_1,\ldots,u_r;q)$-modules
can be constructed with Theorems 2.8 and 3.8.

If $\lambda/\mu$ is a skew shape define
$$NW(\lambda/\mu) = \{\hbox{northwest corner boxes of
$\lambda/\mu$}\},$$
so that $NW(\lambda/\mu)$ is the set of boxes $b\in
\lambda/\mu$ such that there is no box of $\lambda/\mu$ immediately above 
or
immediately to the left of $b$.

\thm
Fix $u_1,\ldots,u_r\in \CC^*$ and let $(c,\lambda/\mu)$ be a placed
skew shape with $n$ boxes.  If
$$\{ q^{2c(b)} \ |\ b\in NW(\lambda/\mu)\} \subseteq \{u_1,\ldots, u_r\}$$
then the $H_{\infty,1,n}$-module $\tilde H^{(c,\lambda/\mu)}$ is a simple
$H_{r,1,n}(u_1,\ldots,u_r;q)$-module (via Theorem 2.8).
\pf
Let $(c,\lambda/\mu)$ be a placed skew shape with $n$ boxes and let
$\tilde H^{(c,\lambda/\mu)}$ be the simple
$H_{\infty,1,n}$-module of Theorem 3.8.  Let
$\rho^{(c,\lambda/\mu)}\colon H_{\infty,1,n}
\to \End(\tilde H^{(c,\lambda/\mu)})$
be the representation corresponding to $\tilde H^{(c,\lambda/\mu)}$.  By 
the
formulas in Theorem 3.8 the matrix $\rho(X^{\varepsilon_1})$ is diagonal 
with
eigenvalues $q^{2c(L(1))}$, for $L\in \cF^{(c,\lambda/\mu)}$, where
$\cF^{(c,\lambda/\mu)}$ is the set of standard tableaux of shape
$\lambda/\mu$.

The boxes $L(1)$, $L\in \cF^{(c,\lambda/\mu)}$ are exactly the northwest
corner boxes of $\lambda/\mu$ and so the minimal polynomial of
$\rho(X^{\varepsilon_1})$ is
$$p(t) = \prod_{b\in NW(\lambda/\mu)} (t-q^{2c(b)}).$$
Thus the condition
$\{ q^{2c(b)} \ |\ b\in NW(\lambda/\mu)\} \subseteq \{u_1,\ldots, u_r\}$
is exactly what is needed for $\tilde H^{(c,\lambda/\mu)}$ to be
an $H_{r,1,n}(u_1,\ldots,u_r;q)$-module via Theorem 2.8.
\endpf

\thm  If the cyclotomic Hecke algebra
$H_{r,1,n}(u_1,u_2,\ldots,u_r;q)$ is semisimple,
then its simple modules are the modules
$\tilde H^{(c,\lambda)}$ constructed in Theorem 3.17,
where
$\lambda=(\lambda^{(1)},\ldots,\lambda^{(r)})$
is an $r$-tuple of partitions with a total of $n$ boxes and
$c$ is the content function determined by
$$q^{2c(b)} = u_i,\quad\hbox{if $b$ is the northwest
corner box of $\lambda^{(i)}$.}$$
\pf
The equations $q^{2c(b)}=u_i$, for $b\in NW(\lambda^{(i)})$,
determine the values $c(b)$ for all boxes $b\in \lambda$
and thus the pair $(c,\lambda)$ defines a placed skew shape.
By a Theorem of Ariki [Ar1],
$H_{r,1,n}(u_1,\ldots,u_r;q)$ is semisimple if and only if
$$[n]_q!\ne 0
\qquad\hbox{and}\qquad u_iu_j^{-1}\not\in\{1, q^2,\ldots,q^{2n}\},$$
where $[n]_q! = [n]_q[n-1]_q\cdots [2]_q[1]_q$ and
$[k]_q = (q^k-q^{-k})/(q-q^{-1})$.
These conditions guarantee that $(c,\lambda)$ is a
placed skew shape and that the $\tilde H_{\infty,1,n}$-module
$\tilde H^{(c,\lambda)}$ defined in Theorem 3.8 is
well defined and irreducible.  The reduction in Theorem 2.8
makes $\tilde H^{(c,\lambda)}$ into a simple
$H_{r,1,n}(u_1,\ldots,u_r;q)$-module and a count of standard
tableaux (using binomial coefficients and the classical identity for the
symmetric
group case) shows that
$$\sum_{\lambda =  (\lambda^{(1)},\ldots,\lambda^{(r)})}
\dim(\tilde H^{(c,\lambda)})^2 = r^nn! = \dim(H_{r,1,n}),$$
where the sum is over all $r$-tuples $\lambda =
(\lambda^{(1)},\ldots,\lambda^{(r)})$ with a total of $n$ boxes.
Thus the $\tilde H^{(c,\lambda)}$ are a complete set of
simple $H_{r,1,n}(u_1,\ldots,u_r;q)$-modules.
\endpf

Theorem 3.18 demonstrates that the construction of simple modules for
$H_{r,1,n}(u_1,\ldots, u_r;q)$ by Ariki and Koike [AK, Theorem 3.7],
for $H_{2,1,n}(p,-p^{-1};q)$ by Hoefsmit [Hf],
for $H_{1,1,n}(1;q)$ by Hoefsmit [Hf] and Wenzl [Wz] (independently),
and for $H_{1,1,n}(1;1)=\CC S_n$ and $H_{2,1,n}(1,1;1)=\CC WB_n$
by Young [Yg1-2], are all special cases of Theorem 3.8.

\subsection Jucys-Murphy elements in cyclotomic Hecke algebras.

The following result is well known, but we give a new
proof which shows that the cyclotomic Hecke algebra analogues of the
Jucys-Murphy elements which have appeared in the literature (see [BMM],
[Ra1], [DJM] and the references there) come naturally from the affine Hecke
algebra $H_{\infty,1,n}$.

\cor
Let $H_{r,1,n}(u_1,\ldots,u_r;q)$ and $H_{r,p,n}(x_0,\ldots,x_{d-1};q)$
be the cyclotomic Hecke algebras defined in (1.1) and (1.2).
\item{(a)}  The elements
$$M_i =  T_i\cdots T_2T_1T_2\cdots T_i\,,\qquad 1\le i\le n,$$
generate a commutative subalgebra of $H_{r,1,n}(u_1,\ldots,u_r;q)$.
\item{(b)}  If $H_{r,1,n}(u_1,\ldots,u_r;q)$ is semisimple then
every simple $H_{r,1,n}(u_1,\ldots,u_r;q)$-module
has a basis of simultaneous eigenvectors of the elements
$M_i$.
\item{(c)}
The elements
$$M_1^p=a_0,\qquad\hbox{and}\qquad
M_iM_1^{-1} = a_i\cdots a_3a_2a_1a_3\cdots a_i\,, \qquad
2\le i\le n,$$
generate a commutative subalgebra of $H_{r,p,n}(x_0,\ldots,x_{d-1};q)$.
\item{(d)}  If $H_{r,p,n}(x_0,\ldots,x_{d-1};q)$ is semisimple then
every simple $H_{r,p,n}(x_0,\ldots,x_{d-1};q)$-module
has a basis of simultaneous eigenvectors of the elements
$M_1^p$ and $M_iM_1^{-1}$, $2\le i\le n$.
\pf
(a)  The elements $X^{\varepsilon_i}$, $1\le i\le n$, generate
the subalgebra $\CC[X]\subseteq H_{\infty,1,n}$.  Inductive
use of the relation (1.4) shows that
$M_i=\Phi(X^{\varepsilon_i})$, where $\Phi\colon H_{\infty,1,n}
\to H_{r,1,n}(u_1,\ldots, u_r;q)$ is the homomorphism in Proposition 2.5.
Thus the subalgebra of $H_{r,1,n}(u_1,\ldots, u_r;q)$ generated
by the element $M_i$ is the image of $\CC[X]$ under the homomorphism
$\Phi$.

\smallskip\noindent
(b) is an immediate consequence of Theorem 3.18, the construction
described in
Theorem 3.8, and the fact that $M_i = \Phi(X^{\varepsilon_i})$.

\smallskip\noindent
(c)  The elements $X^{p\varepsilon_1}$ and
$X^{\varepsilon_i-\varepsilon_1}$, $2\le i\le n$, generate the subalgebra
$\CC[X^{L_p}]\subseteq H_{\infty,p,n}$ and the images of these
elements under the homomorphism $\Phi_p\colon H_{\infty,p,n}\to
H_{r,p,n}(x_0,\ldots, x_{d-1};q)$ are the elements $M_1^p$ and
$M_iM_1^{-1}$.
\smallskip\noindent
The proof of part (d) uses the construction
described in
Theorem 3.15 and is analogous to the proof of part~(b).
\endpf

\subsection  The center of $H_{r,1,n}(u_1,\ldots,u_r;q)$.

It is an immediate consequence of (1.12) and the proof of Corollary 3.20
that
$$\CC[M_1,\ldots,M_n]^{S_n}\subseteq
Z(H_{r,1,n}(u_1,\ldots,u_r;q)),$$
where $M_i=T_iT_{i-1}\cdots T_2T_1T_2\cdots T_{i-1}T_i$.
The following proposition shows that
this inclusion is an equality.

\prop If the cyclotomic Hecke algebra
$H_{r,1,n}(u_1,\ldots,u_r;q)$ is semisimple then its center
$$
Z(H_{r,1,n}(u_1,\ldots,u_r;q))=
\CC[M_1,\ldots,M_n]^{S_n},
$$
where $M_i=T_iT_{i-1}\cdots T_2T_1T_2\cdots T_{i-1}T_i$
and $\CC[M_1,\ldots,M_n]^{S_n}$ is the ring of symmetric polynomials in
$M_1,\ldots, M_n$.
\pf
By (1.12) $Z(H_{\infty,1,n}) = \CC[X]^{S_n}$.
Thus, since $M_i=\Phi(X^{\varepsilon_i})$,
where $\Phi\colon H_{\infty,1,n}\to H_{r,1,n}(u_1,\ldots,u_r;q)$
is the surjective homomorphism of Proposition 2.5, it follows that
$$\CC[M_1,\ldots,M_n]^{S_n}=\Phi(\CC[X]^{S_n})
=\Phi(Z(H_{\infty,1,n}))\subseteq Z(H_{r,1,n}(u_1,\ldots,u_r;q)).$$

For the reverse inclusion we need to show that the
action of the elements $\CC[M_1,\ldots, M_n]^{S_n}$
distinguishes the simple $H_{r,1,n}(u_1,\ldots,u_r;q)$-modules.

Let $\lambda=(\lambda^{(1)},\ldots,\lambda^{(r)})$ be an $r$-tuple of
partitions with
a total of $n$ boxes and let $\tilde H^{(c,\lambda)}$
be the corresponding simple $H_{r,1,n}(u_1,\ldots,u_r;q)$-module
as constructed by Theorem 3.18 (and Theorem 3.8).
If $L$ is a standard tableau of shape
$\lambda$ then $M_1,\ldots, M_n$ act on $v_L$ by the multiset of
values $(q^{2c(L(1))}, \ldots, q^{2c(L(n))})$.
The elementary symmetric functions $e_i(M_1,\ldots, M_n)$
act on $\tilde H^{(c,\lambda)}$ by the values
$a_i=e_i(q^{2c(L(1))}, \ldots, q^{2c(L(n))})$.  Note that $a_i$
does not depend on the choice of the standard tableau $L$
(since $e_i(M_1,\ldots, M_n)\in Z(H_{r,1,n})$).
We show that the simple module $\tilde H^{(c,\lambda)}$
is determined by the values $a_1,\ldots,a_n$.  This
shows that the simple modules are distinguished by the elements
of $\CC[M_1,\ldots, M_n]^{S_n}$.

Let us explain how
the values $a_1,\ldots,a_n$ determine the placed skew shape $(c,\lambda)$.
There is a unique (unordered) multiset of values $b_1,\ldots, b_n$
such that $e_i(b_1,\ldots, b_n)=a_i$ for all $1\le i\le n$.  The
$b_i$ are determined by the equation
$$(t-b_1)\cdots (t-b_n) = t^n-a_1t^{n-1}+a_2t^{n-2}
-a_3t^{n-3}+\cdots\pm a_n.$$
  From the previous paragraph, the multiset
$\{b_1,\ldots,b_n\}$ must be the same as the multiset

\noindent
$\{q^{2c(L(1))}, \ldots, q^{2c(L(n))}\}$.
In this way the multiset $S=\{c(L(1)),\ldots, c(L(n))\}$
is determined by $a_1,\ldots, a_n$.

The multiset $S$ is a disjoint union
$S=S_1\sqcup\cdots\sqcup S_r$ of multisets such that each $S_i$
is in a single $\ZZ$-coset of $\RR+i[0,2\pi/\ln(q^2))i$ (see (3.4)).
These $\ZZ$ cosets are determined by the values $u_1,\ldots, u_r$
and thus there is a one-to-one correspondence between the
$u_i$ and the multisets $S_i$.
Then $|\lambda^{(i)}|={\rm Card}(S_i)$,
the nonempty diagonals of $\lambda^{(i)}$
are determined by the values in $S_i$,
and the lengths of the diagonals
of $\lambda^{(i)}$ are the multiplicities of the values
of the elements of $S_i$.  This information completely determines
$\lambda^{(i)}$ for each $i$.  Thus $(c,\lambda)$ is determined
by the values $a_1,\ldots, a_n$.
\endpf

\remark In the language of affine Hecke algebra representations
(see [Ra2]) the proof of Proposition 3.22 shows that (when
$H_{r,1,n}(u_1,\ldots,u_r;q)$ is semisimple)
the simple $H_{r,1,n}(u_1,\ldots,u_r;q)$-modules $\tilde H^{(c,\lambda)}$
all have different central characters (as $H_{\infty,1,n}$ modules).
\endprop

\remark
The elements $M_1^p$ and $M_iM_1^{-1}$ from Corollary
3.20 cannot be used to obtain a direct analogue of Proposition 3.22 for
$H_{r,p,n}(x_0,\ldots,x_{d-1};q)$.  This is because all of the
$H_{\infty,p,n}$-modules $V_j$ appearing in the decomposition (3.14)
will have the same central character.
However, {\it when $n$ is odd}, the natural analogue of Proposition 3.22 
does
hold for Iwahori-Hecke algebras of type $D_n$, $HD_n(q)=H_{2,2,n}(1;q)$.
In that case, every simple $H_{2,1,n}$-module has trivial inertia
group and the decomposition in (3.14) has only one summand.

\section 4. Affine Hecke algebras of general type

Let $R$ be a reduced root system and let $R^\vee$ be the root system formed
by the coroots
$\alpha^\vee = {2\alpha/ \langle \alpha,\alpha\rangle}$,
for $\alpha\in R$.
Let $W$ be the Weyl group of $R$
and fix a system of positive roots $R^+$ in $R$. Let $\{\alpha_1,
\ldots,\alpha_n\}$ be the corresponding simple roots and let
$s_1,\ldots, s_n$ be the corresponding simple reflections in $W$.
The fundamental weights are defined by the equations
$\langle \omega_i,\alpha_j^\vee\rangle = \delta_{ij}$ and the lattices
$$
P=\sum_{i=1}^n \ZZ\omega_i
\quad\hbox{and}\quad
Q=\sum_{i=1}^n \ZZ\alpha_i,$$
are the weight lattice and the root lattice, respectively.
The Dynkin diagrams and the corresponding
extended Dynkin diagrams are given in Figure 3.  If $\Gamma$ is a
Dynkin diagram or extended Dynkin diagram define
$$m_{ij}= \cases{
2, &if
\beginpicture
\setcoordinatesystem units <1cm,1cm>       % sets scale
\setplotarea x from 2.8 to 4.2, y from 1.8 to 2 % sets plot size up
{\scriptsize
\multiput {$\circ$} at 3  1.9 *1 1 0 /      % puts nodes in
\put {$^i$}     at 3 2.1   %
\put {$^j$}     at 4 2.1   % label nodes with roots above
\linethickness=0.5pt                      % sets line thickness
}
\endpicture
,
\cr \cr
3, &if
\beginpicture
\setcoordinatesystem units <1cm,1cm>       % sets scale
\setplotarea x from 2.8 to 4.2, y from 1.8 to 2 % sets plot size up
{\scriptsize
\multiput {$\circ$} at 3  1.9 *1 1 0 /      % puts nodes in
\put {$^i$}     at 3 2.1   %
\put {$^j$}     at 4 2.1   % label nodes with roots above
\linethickness=0.5pt                      % sets line thickness
\putrule from 3.05 1.9 to 3.95 1.9         %   puts solid lines between 
nodes
}
\endpicture
,
\cr}
\qquad\hbox{and}\qquad
m_{ij} = \cases{
4, &if
\beginpicture
\setcoordinatesystem units <1cm,1cm>       % sets scale
\setplotarea x from 2.8 to 4.2, y from 1.8 to 2 % sets plot size up
{\scriptsize
\multiput {$\circ$} at 3  1.9 *1 1 0 /      % puts nodes in
\put {$^i$}     at 3 2.1   %
\put {$^j$}     at 4 2.1   % label nodes with roots above
\linethickness=0.5pt                      % sets line thickness
\putrule from 3.03 1.935 to 3.97 1.935  %
\putrule from 3.03 1.865 to 3.97 1.865  %
}
%\put{Undirected Dynkin diagram of type $G_2$}[b] at 3.5 1.25
\endpicture
,
  \cr \cr
6, &if
\beginpicture
\setcoordinatesystem units <1cm,1cm>       % sets scale
\setplotarea x from 2.8 to 4.2, y from 1.8 to 2 % sets plot size up
{\scriptsize
\multiput {$\circ$} at 3  1.9 *1 1 0 /      % puts nodes in
\put {$^i$}     at 3 2.1   %
\put {$^j$}     at 4 2.1   % label nodes with roots above
\linethickness=0.5pt                      % sets line thickness
\putrule from 3.03 1.945 to 3.97 1.945  %
\putrule from 3.03 1.855 to 3.97 1.855  %
\putrule from 3.05 1.9 to 3.95 1.9         %   puts solid lines between 
nodes
}
%\put{Undirected Dynkin diagram of type $G_2$}[b] at 3.5 1.25
\endpicture
.
\cr}$$

\subsection Affine Weyl groups.

The extended affine Weyl group is the group
$$
\tilde W = W\ltimes P =\{ wt_\lambda\ |\ w\in W,\ \lambda\in P\},
\quad\hbox{where}\quad
wt_\lambda = t_{w\lambda}w,
$$
for $w\in W$ and $\lambda\in P$
where $t_\lambda$ corresponds to translation by $\lambda\in P$.
Define $s_0\in \tilde W$ by the equation
$$s_0s_{\phi^\vee} = t_{\phi},
\qquad\hbox{where $\phi^\vee$ is the highest root of $R^\vee$,}
\formula$$
see [Bou Ch. IV \S 1 no. 2.1].
The subgroup
$W_{\rm aff} = W\ltimes Q$ of $\tilde W$
is presented by generators
$s_0,s_1,\ldots,s_n$ and relations
$$
s_i^2 =1, \quad 0\le i\le n,
\qquad\hbox{and}\qquad
\underbrace{s_is_j\cdots}_{m_{ij}} = \underbrace{s_js_i\cdots}_{m_{ij}},
\quad i\ne j,
$$
where the $m_{ij}$ are determined from the extended Dynkin diagram
of the root system $R^\vee$.   Define
$$\Omega=\{g_i\ |\ \omega_i\hbox{\ is minuscule}\},
\qquad\hbox{where}\quad g_iw_0w_{0,i}= t_{\omega_i},
\formula$$
$w_0$ is the longest element of $W$ and $w_{0,i}$ is the longest element
of the group $\langle s_j\ |\ 1\le j\le n, j\ne i\rangle$, see [Bou Ch. IV 
\S 2
Prop. 6]. Then $\Omega\cong P/Q$ and each element $g\in \Omega$ 
corresponds to
an automorphism of the extended Dynkin diagram of $R^\vee$, in the
sense that
$$\hbox{if $g\in \Omega$\quad then\quad $gs_ig^{-1} = s_{\sigma(i)}$,}
\formula$$
where $\sigma$ is the permutation of the nodes determined by the
automorphism.  Equation (4.4) means that $\tilde W = W_{\rm 
aff}\rtimes\Omega$.
The usual length function on the Coxeter group
$W_{\rm aff}$ is extended to the group $\tilde W$ by
$$
\ell(wg) = \ell(w), \qquad\hbox{for $w\in W_{\rm aff}$ and $g\in \Omega$.}
$$

Let $L$ be a lattice such that $Q\subseteq L\subseteq P$.
View $L/Q$ as a subgroup of $\Omega\cong P/Q$ and let
$$
\tilde W_L = W\ltimes L = W_{\rm aff}\rtimes (L/Q).
$$
Then $W_{\rm aff} = \tilde W_Q$, $\tilde W= \tilde W_P$, and
$\tilde W_L$ is a subgroup of $\tilde W$.

\subsection Affine braid groups.

Let $L$ be a lattice such that $Q\subseteq L\subseteq P$.
The affine braid group $\tilde\cB_L$ is the group
given by generators
$T_w$, $w\in \tilde W_L$, and relations
$$
T_wT_{w'} = T_{ww'},\qquad\hbox{if $\ell(ww')=\ell(w)+\ell(w')$.}
$$
Let $\tilde \cB_{\rm aff}=\tilde \cB_Q$.
View $L/Q$ as a subgroup of $\Omega\cong P/Q$. Then
$\cB_L = \cB_{\rm aff}\rtimes L/Q$ is presented by generators
$T_i=T_{s_i}$, $0\le i\le n$, and relations
$$\underbrace{T_iT_j\cdots}_{m_{ij}} = \underbrace{T_jT_i\cdots}_{m_{ij}},
\qquad\hbox{and}
\qquad
gT_ig^{-1}=T_{\sigma(i)},\quad\hbox{for $g\in\Omega$},
\formula$$
where $\sigma$ is as in (4.4), and the $m_{ij}$ are specified by
the extended Dynkin diagram of $R^\vee$.

Let $P^+ = \sum_{i=1}^n \ZZ_{\ge 0}\omega_i$ be the dominant weights in 
$P$.
Define elements $X^\lambda$, $\lambda\in P$ by
$$
X^\lambda = T_{t_\lambda},\quad\hbox{if $\lambda\in P^+$, and}
\qquad
X^\lambda = X^\mu(X^\nu)^{-1},
\quad\hbox{if
$\lambda=\mu-\nu$ with $\mu,\nu\in P^+$.}
\formula$$
By [Mac3, 3.4] and [Lu], the $X^\lambda$ are well defined and
do not depend on the choice $\lambda=\mu-\nu$, and
$$X^\lambda X^\mu = X^\mu X^\lambda = X^{\lambda+\mu},
\quad \hbox{for $\lambda,\mu\in P$.}
\formula$$
Then $X^\lambda\in \tilde \cB_L$ if and only if $\lambda\in L$.

\subsection Affine Hecke algebras.

Fix $q\in\CC^*$.
The affine Hecke algebra $\widetilde H_L$
is the quotient of the group algebra $\CC\tilde\cB_L$ by the
relations
$$(T_i-q)(T_i+q^{-1})  = 0,\qquad\hbox{for $0\le i\le n$.}\formula$$
In $\tilde H_L$ (see [Mac, 4.2]),
$$
X^\lambda T_i = T_iX^{s_i\lambda}
+ (q-q^{-1}){X^\lambda-X^{s_i\lambda}\over 1-X^{-\alpha_i}},
\qquad\hbox{for $\lambda\in L$, $1\leq i\leq n$.}
\formula$$
The Iwahori-Hecke algebra $H$ is the
subalgebra of $\tilde H$ generated by $T_1,\ldots, T_n$.

\medskip
To our knowledge, the following theorem is due to Bernstein and
Zelevinsky in type A, and to Bernstein in general type (unpublished).
Lusztig has given an exposition in [Lu].
We give a new proof which we believe is more elementary and more direct.

\thm  (Bernstein, Zelevinsky, Lusztig [Lu])
Let $L$ be a lattice such that $Q\subseteq L\subseteq P$, where $Q$ is the
root lattice and $P$ is the weight lattice of the root system $R$.
Let $\tilde H = \tilde H_L$ be the affine Hecke algebra corresponding
to $L$ and let $\CC[X] ={\rm span}\{ X^\lambda\ |\ \lambda\in L\}$.
Let $W$ be the Weyl group of $R$.  Then the center of $\tilde H$ is
$$
Z(\tilde H)= \CC[X]^{W}
=\left\{ f\in\CC[X]\ |\ w f = f
\hbox{ for every } w\in W\right\}.
$$
\vskip-.1in
\pf
Assume
$$
z=\sum_{\lambda\in L,w\in W} c_{\lambda,w}X^\lambda T_w\in Z(\tilde H).
$$
Let $m\in W$ be maximal in Bruhat order subject to $c_{\gamma,m}\neq0$ for
some $\gamma\in L$.
If $m\ne 1$ there exists a dominant $\mu\in L$ such that
$c_{\gamma+\mu-m\mu,m}=0$ (otherwise  $c_{\gamma+\mu-m\mu,m}\neq 0$ for 
every
dominant $\mu\in L$, which is impossible since $z$ is a finite linear
combination of $X^\lambda T_w$). Since $z\in Z(\tilde H)$
we have
$$
z  =   X^{-\mu}zX^\mu   =
\sum_{\lambda\in L,w\in W} c_{\lambda,w} X^{\lambda-\mu} T_w X^\mu.
$$
Repeated use of the relation (4.11) yields
$$
T_wX^\mu=\sum_{\nu\in L,v\in W} d_{\nu,v}X^\nu T_v
$$
where $d_{\nu,v}$ are constants such that
$d_{w\mu,w}=1$, $d_{\nu,w}=0$ for
$\nu\ne w\mu$, and $d_{\nu,v}=0$ unless $v\le w$.
So
$$
z  = \sum_{\lambda\in L,w\in W} c_{\lambda,w}X^\lambda T_w
= \sum_{\lambda\in L,w\in W}\sum_{\nu\in L,v\in W}
c_{\lambda,w}d_{\nu,v} X^{\lambda-\mu+\nu} T_v
$$
and comparing the coefficients of $X^\gamma T_m$ gives
$
c_{\gamma,m}=c_{\gamma+\mu-m\mu,m} d_{m\mu,m}.
$
Since $c_{\gamma+\mu-m\mu,m}=0$ it follows that
$c_{\gamma,m}=0$, which is a contradiction.  Hence
$z=\sum_{\lambda\in L} c_\lambda X^\lambda\in \CC[X]$.

The relation (4.11) gives
$$
zT_i=T_iz=(s_iz)T_i+(q-q^{-1})z'
$$
where $z'\in\CC[X]$. Comparing coefficients of
$X^\lambda$  on both sides yields $z' = 0$. Hence
$zT_i=(s_iz)T_i$, and therefore $z=s_iz$ for $1\leq i\leq n$.  So
$z\in \CC[X]^W$.
\endpf

\subsection Deducing the $\tilde H_L$ representation theory from $\tilde 
H_P$.

Although we have not taken this point of view in our presentation,
  the affine Hecke algebras defined above are naturally
associated to a reductive algebraic group $G$ over $\CC$ [KL] or a
$p$-adic Chevalley group [IM].  In this formulation,
the lattice $L$ is determined by the group of characters of
the maximal torus of $G$.
It is often convenient to work only with the adjoint version or
only with the simply connected version of the group $G$ and
therefore it seems desirable to be able to derive the
representation theory of the affine Hecke algebras $\tilde H_L$
from the representation theory of the affine Hecke algebra $\tilde H_P$.
The following theorem shows that this can be done in a simple
way by using the extension of Clifford theory in
the Appendix.  In particular, Theorem A.13, can be used to construct
all of the simple $\tilde H_L$-modules from the simple $\tilde 
H_P$-modules.

\thm
Let $L$ be a lattice such that $Q\subseteq L\subseteq P$, where $Q$ is the
root lattice and $P$ is the weight lattice of the root system $R$.
Let $\tilde H = \tilde H_L$ be the affine Hecke algebra corresponding
to $L$.
Then there is an action of a finite group $K$ on $\tilde H_P$,
acting by automorphisms, such that
$$\tilde H_L = (\tilde H_P)^K,$$
is the subalgebra of fixed points
under the action
of the group $K$.
\pf
There are two cases to consider, depending on whether the
group $\Omega\cong P/Q$ is cyclic or not.
$$
\matrix{
\hbox{Type} & A_{n-1} & B_n & C_n & D_{2n-1} & D_{2n} &
E_6 & E_7 & E_8 & F_4 & G_2 \cr
\Omega & \ZZ/n\ZZ & \ZZ/2\ZZ & \ZZ/2\ZZ & \ZZ/4\ZZ & 
\ZZ/2\ZZ\times\ZZ/2\ZZ&
\ZZ/3\ZZ& \ZZ/2\ZZ &  1 & 1 & 1 \cr
}$$
In each case we construct the group $K$ and its action on
$\tilde H_P$ explicitly.  This is necessary for the effective application
of Theorem A.13 on examples.

\medskip\noindent
{\it Case 1.}  If $\Omega$ is a cyclic group
$\Omega$ then the subgroup $L/Q$ is a cyclic subgroup. Suppose
$$\Omega = \{1,g,\ldots, g^{r-1}\}\qquad\hbox{and}\qquad
L/Q = \{1, g^d, \ldots, g^{d(p-1)}\},$$
where $pd=r$.
Let $\zeta$ be a primitive $p$th root of unity and define an
automorphism
$$\matrix{
\sigma\, \colon &\tilde H_P &\longrightarrow &\tilde H_P\cr
&g &\longmapsto &\zeta g\, ,\cr
&T_i &\longmapsto &T_i\, , &\qquad &0\le i\le n.\cr}$$
The map $\sigma$ is an algebra isomorphism since it preserves the relations
in (4.6) and (4.10).  Furthermore, $\sigma$ gives rise to a $\ZZ/p\ZZ$
action on
$\tilde H_P$ and
$$\tilde H_L = (\tilde H_P)^{\ZZ/p\ZZ}.\formula$$

\smallskip\noindent
{\it Case 2.}  If the root system $R^\vee$ is of type $D_n$, $n$ even,
then $\Omega\cong \ZZ/2\ZZ\times\ZZ/2\ZZ$ and the subgroups of $\Omega$
correspond
to the intermediate lattices $Q\subseteq L\subseteq P$.
Suppose
$$\Omega=\{1,g_1,g_2,g_1g_2\ |\ g_1^2=g_2^2=1, g_1g_2=g_2g_1\}.$$
and define automorphisms of $\tilde H_P$ by
$$\matrix{
\sigma_1\, \colon &\tilde H_P &\longrightarrow &\tilde H_P\cr
&g_1 &\longmapsto &-g_1\, ,\cr
&g_2 &\longmapsto &g_2\, ,\cr
&T_i &\longmapsto &T_i\, , \cr}
\qquad\hbox{and}\qquad
\matrix{
\sigma_2\, \colon &\tilde H_P &\longrightarrow &\tilde H_P\cr
&g_1 &\longmapsto &g_1\, ,\cr
&g_2 &\longmapsto &-g_2\, ,\cr
&T_i &\longmapsto &T_i\, .\cr}
$$
Then
$$\tilde H_{L_1} = (\tilde H_P)^{\sigma_1},
\qquad
\tilde H_{L_2} = (\tilde H_P)^{\sigma_2},
\qquad\hbox{and}\qquad
\tilde H_Q = (\tilde H_P)^{\langle \sigma_1,\sigma_2\rangle},
\formula$$
where $L_1$ and $L_2$ are the two intermediate lattices strictly between
$Q$ and
$P$.
\endpf

\section 5. Where does the homomorphism $\Phi$ come from?

The homomorphism
$\Phi\colon H_{\infty,1,n}\to H_{r,1,n}$ of
Proposition 2.5 is a powerful tool for transporting results about the
affine Hecke algebra of type $A$ to the cyclotomic Hecke algebras.
In this section we show how this homomorphism arises naturally,
from a folding of the Dynkin diagram of $\tilde B_n$,
and we give some generalizations of the homomorphism
$\Phi$ to other types.

\bigskip\noindent
{\bf Example 1. Type $C_n$.}
The root system $R$ of type $C_n$ can be realized by
$$R=\{\pm 2\varepsilon_i, \pm(\varepsilon_j-\varepsilon_i)
\ |\ 1\le i,j\le n\},$$
where $\varepsilon_i$ are an orthonormal basis of $\RR^n$.  The simple
roots and the fundamental weights are given by
$$\eqalign{
\alpha_1 &= 2\varepsilon_1, \qquad
\alpha_i = \varepsilon_i-\varepsilon_{i-1},
\qquad 2\le i\le n,\cr \cr
\omega_i &= \varepsilon_n + \varepsilon_{n-1}+\cdots+\varepsilon_i,
\qquad 1\le i\le n. \cr
}$$
If $\phi^\vee$ is the highest root of $R^\vee$ then
$$\phi = \varepsilon_n+\varepsilon_{n-1},
\qquad\hbox{and}\qquad
s_{\phi^\vee} = (s_{n-1}\cdots s_2s_1s_2\cdots s_{n-1})s_n
(s_{n-1}\cdots s_2s_1s_2\cdots s_{n-1}).
\formula$$
Then $\omega_n = \varepsilon_n$ is the only miniscule weight,
$$\eqalign{w_0 &= (1,-1)(2,-2)\cdots (n,-n),\cr
w_{0,n} &= (1,-1)(2,-2)\cdots (n-1,-(n-1)),\qquad\hbox{and}\cr
w_0w_{0,n} &= (n,-n) = s_n\cdots s_2s_1s_2\cdots s_n.\cr
}
\formula$$
Thus, from (4.2), (4.3) and (4.6), $\Omega = \{1, g_n\}\cong \ZZ/2\ZZ$,
$$\eqalign{
g_n &= X^{\varepsilon_n}
T_n^{-1}\cdots T_2^{-1}T_1^{-1}T_2^{-1} \cdots T_n^{-1}, \cr
T_0 &= X^{\varepsilon_n+\varepsilon_{n-1}}
(T_n^{-1}\cdots T_2^{-1}T_1^{-1}T_2^{-1}\cdots T_n^{-1})
T_n^{-1}
(T_n^{-1}\cdots T_2^{-1}T_1^{-1}T_2^{-1}\cdots T_n^{-1}),\cr}
\formula$$
and
$$g_nT_0g_n^{-1}=T_n, \qquad\hbox{and}\qquad
g_nT_ng_n^{-1} = T_0.\formula$$

The braid group $\tilde \cB_P(C_n)$ is generated by $T_0,T_1,\ldots, T_n$ 
and
$g_n$ which satisfy relations in (4.6), where the $m_{ij}$ are given by
the extended Dynkin diagram $\tilde B_n$, see Figure 3.
The braid group $\cB(C_n)$ is the subgroup
generated by $T_1,\ldots, T_n$.  These elements
satisfy the relations in (4.6), where the $m_{ij}$ are given by the Dynkin
diagram
$C_n$. A straightforward check verifies that the map defined by
$$
\beginpicture
\setcoordinatesystem units <1cm,1cm>        % sets scale
%********************************************************************
% Dynkin diagram of type extended B_n
%********************************************************************
\setplotarea x from 2.8 to 9.7, y from -0.8 to 0.8  % sets plot size up
{\scriptsize
\multiput {$\circ$} at 9 -0.5 *1 0 1 /      %
\multiput {$\circ$} at 7 0    *1 1 0 /      %
\multiput {$\circ$} at 3 0    *2 1 0 /      %  puts nodes in
\put {$^0$}     at 9 -0.9   %
\put {$^n$}     at 9 0.7   %
\put {$^{n-1}$} at 7.9 0.2 %
\put {$^{n-2}$} at 7 0.2   %
\put {$^3$}     at 5 0.2   %
\put {$^2$}     at 4 0.2   %
\put {$^1$}     at 3 0.2   % label nodes with roots below
\setplotsymbol ({\rm .})
\setquadratic \plot 9.2 0.5                      % draws automorphism
                     9.5 0                        %
                     9.2 -0.5 /                    %
\arrow <5pt> [.2,.67] from 9.3 0.4    to 9.2 0.5  %
\arrow <5pt> [.2,.67] from 9.3 -0.4  to 9.2 -0.5  %
\linethickness=0.75pt                          % sets line thickness
\putrule from 3.03 0.045 to 3.97 0.045       % puts solid lines between 
nodes
\putrule from 3.03 -0.045 to 3.97 -0.045       %
\putrule from 4.05 0   to 4.95 0               %
\putrule from 7.05 0   to 7.95 0               %
\setlinear \plot 8.95 -0.5   8.05 -0.05 /     %
\setlinear \plot 8.95 0.5   8.05 0.05 /     %
\setdashes <2mm,1mm>          %
\putrule from 5.05 0   to 6.95 0    % draws dotted lines between nodes
}
% end of dynkin diagram of type extended B_n
\endpicture
\qquad
\longrightarrow
\qquad
%********************************************************************
% Dynkin diagram of type B_n
%********************************************************************
\beginpicture
\setcoordinatesystem units <1cm,1cm>       % sets scale
\setplotarea x from 2.8 to 8.2, y from -0.8 to 0.8 % sets plot size up
{\scriptsize
\multiput {$\circ$} at 3   0 *2 1 0 /      %
\multiput {$\circ$} at 7   0 *1 1 0 /      %  puts nodes in
\put {$^1$}     at 3 0.2   %
\put {$^2$}     at 4 0.2   %
\put {$^3$}     at 5 0.2   %
\put {$^{n-1}$}   at 7 0.2   %
\put {$^n$}     at 8 0.2   % label nodes with roots above
\linethickness=0.75pt                      % sets line thickness
\putrule from 3.03 0.045 to 3.97 0.045  %
\putrule from 3.03 -0.045 to 3.97 -0.045  %
\putrule from 4.05 0 to 4.95 0         %   puts solid lines between nodes
\putrule from 7.05 0 to 7.95 0         %
\setdashes <2mm,1mm>            %
\putrule from 5.05 0 to 6.95 0  % draws dotted lines between nodes
}
%\put{Undirected Dynkin diagram of type $C_n$}[b] at 5.5 1.25
\endpicture
$$
$$\hskip1.2in\matrix{
\Phi_{\tilde CC}\colon &\tilde \cB_P(B_n) &\longrightarrow &\cB(C_n) \cr
&g_n &\longmapsto &1,\cr
&T_0 &\longmapsto &T_n,\cr
&T_i&\longmapsto &T_i, &\qquad &1\le i\le n. \cr
}\formula$$
extends to a well defined surjective group homomorphism.
  From the identity (4.3),
$$\eqalign{
\Phi_{\tilde CC}(X^{\varepsilon_n})
&=\Phi_{\tilde CC}(X^{\omega_n}) \cr
&= \Phi_{\tilde CC}(g_nT_nT_{n-1}\cdots T_2T_1T_2\cdots T_{n-1}T_n)\cr
&=T_nT_{n-1}\cdots T_2T_1T_2\cdots T_{n-1}T_n. \cr}$$
By inductively applying the relation
$X^{\varepsilon_i}=T_iX^{\varepsilon_{i-1}}T_i$ we get
$$\Phi_{\tilde CC}(X^{\varepsilon_i})
= T_iT_{i-1}\cdots T_2T_1T_2\cdots T_{i-1}T_i\,,
\qquad\hbox{for all $1\le i\le n$.\qquad\qed}$$

\bigskip\noindent
{\bf Example 2.  Type $A_{n-1}$.}
Since the weight lattice $P$ for the root system of type $C_n$ is
the  same as the lattice $L$ defined in (1.5)
we have an injective homomorphism
$$\matrix{
\Phi_{\tilde A\tilde C}\colon &\cB_{\infty,1,n} &\longrightarrow
&\tilde \cB_P(C_n)
\cr &T_i&\longmapsto &T_i, &\qquad &2\le i\le n, \cr
&X^{\varepsilon_i} &\longmapsto &X^{\varepsilon_i}.\cr
}$$
The composition of $\Phi_{\tilde A\tilde C}$ and the map
$\Phi_{\tilde CC}$ from (5.5)
is the surjective homomorphism defined by
$$\matrix{
\Phi\colon &\cB_{\infty,1,n} &\longrightarrow &\cB(C_n)\hfill \cr
&T_i&\longmapsto &T_i,\hfill &\qquad &2\le i\le n, \cr
&X^{\varepsilon_i} &\longmapsto
&T_iT_{i-1}\cdots T_2T_1T_2\cdots T_{i-1}T_i,\hfill
&\quad &1\le i\le n.\cr
}$$
In fact, it follows from the defining relations of
$\cB_{\infty,1,n}$ and $\cB(C_n)$ that the map $\Phi$ is an isomorphism!

The cyclotomic Hecke algebras $H_{r,1,n}(u_1,\ldots,u_r;q)$ are quotients
of $\CC \cB(C_n)$ and in this way the
group homomorphism $\Phi$ is the source of the algebra homomorphism
$$\Phi\colon H_{\infty,1,n}\longrightarrow H_{r,1,n}(u_1,\ldots,u_r;q)$$
which was used extensively in Section 3 to relate the representation theory
of the cyclotomic Hecke algebras $H_{r,1,n}(u_1,\ldots,u_r;q)$ to the 
affine
Hecke algebra of type A.  \qquad \qquad\qed

\bigskip\noindent
{\bf Example 3. Type $D_n$.}
Let $R$ be the root system of type $D_n$.  Then
$R^\vee$ is also of type $D_n$ and inspection of the Dynkin diagrams
of types $\tilde D_n$ and $D_n$ yields a surjective algebra homomorphism
defined by
$$
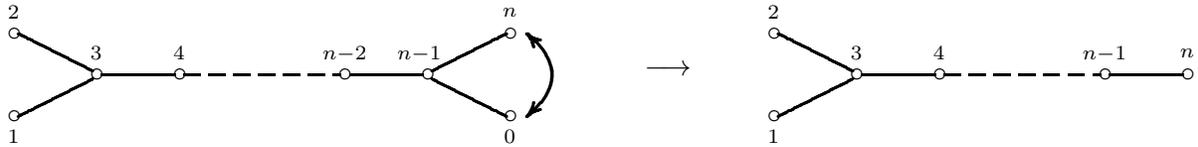
\beginpicture
\setcoordinatesystem units <1cm,1cm>          % sets scale
%********************************************************************
% Dynkin diagram of type extended D_n
%********************************************************************
\setplotarea x from 1.8 to 8.7, y from -0.8 to 0.8 % sets plot size up
{\scriptsize
\multiput {$\circ$} at 2 -0.5 *1 0 1 /      %
\multiput {$\circ$} at 3   0 *1 1 0 /      %
\multiput {$\circ$} at 6   0 *1 1 0 /      %  puts nodes in
\multiput {$\circ$} at 8 -0.5 *1 0 1 /      %
\put {$^0$}   at 8 -0.8   %
\put {$^n$}   at 8  0.7   %
\put {$^{n-1}$} at 6.9 0.2 %
\put {$^{n-2}$} at 6 0.2   %
\put {$^4$}   at 4 0.2   %
\put {$^3$}   at 3 0.2  %
\put {$^2$}   at 2 0.7   %
\put {$^1$}   at 2 -0.8   % label nodes with roots below
\setplotsymbol ({\rm .})
\linethickness=0.75pt                           % sets line thickness
\putrule from 3.05 0 to 3.95 0       % puts solid lines between nodes
\putrule from 6.05 0 to 6.95 0       %
\setlinear \plot 7.05 -0.05  7.95 -0.5 / %
\setlinear \plot 7.05 0.05   7.95 0.5 / %
\setlinear \plot 2.05 -0.5   2.95 -0.05 / %
\setlinear \plot 2.05 0.5   2.95 0.05 / %
\setquadratic \plot 8.2 0.5                      % draws automorphism
                     8.5 0                        %
                     8.2 -0.5 /                    %
\arrow <5pt> [.2,.67] from 8.3 0.4  to 8.2 0.5   %
\arrow <5pt> [.2,.67] from 8.3 -0.4 to 8.2 -0.5   %
\setdashes <2mm,1mm>          %
\putrule from 4.05 0 to 5.95 0  % draws dotted lines between nodes
}
% end of dynkin diagram of type extended D_n
\endpicture
\qquad
\longrightarrow
\qquad
\beginpicture
\setcoordinatesystem units <1cm,1cm>                % sets scale
%********************************************************************
% Dynkin diagram of type D_n
%********************************************************************
\setplotarea x from 2.8 to 8.2, y from -0.8 to 0.8  % sets plot size
{\scriptsize
\multiput {$\circ$} at 8   0 *1 -1 0 /      %
\multiput {$\circ$} at 5   0 *1 -1 0 /      %  puts nodes in
\multiput {$\circ$} at 3 -0.5 *1 0 1 /      %
\put {$^n$}   at 8 0.2   %
\put {$^{n-1}$} at 7 0.2   %
\put {$^4$}   at 5 0.2   %
\put {$^3$}   at 4 0.2   %
\put {$^2$}   at 3 0.7   %
\put {$^1$}   at 3 -0.8  % label nodes with roots below
\setplotsymbol ({\rm .})
\linethickness=0.75pt                           % sets line thickness
\putrule from 7.05 0 to 7.95 0       % puts solid lines between nodes
\putrule from 4.05 0 to 4.95 0       %
\setlinear \plot 3.05 0.5   3.95 0.05 / %
\setlinear \plot 3.05 -0.5   3.95 -0.05 / %
\setdashes <2mm,1mm>          %
\setdashes <2mm,1mm>          %
\putrule from 5.05 0 to 6.95 0  % draws dotted lines between nodes
}
% end of dynkin diagram of type D_n
\endpicture
$$
$$\hskip2in\matrix{
\Phi_{\tilde DD}\colon &\tilde\cB_Q(D_n) &\longrightarrow &\cB(D_n) \cr
&g_n &\longmapsto &1,\cr
&T_0 &\longmapsto &T_n,\cr
&T_i&\longmapsto &T_i, &\qquad &1\le i\le n. &
\qquad\hbox{\qed}\cr
}$$

\medskip
Examples 1, 2, and 3 show that, for types $A$, $B$ and $D$,
there exist surjective homomorphisms from the
affine Hecke algebra to the corresponding Iwahori-Hecke subalgebra.
The following example shows that this is not a general phenomenon:
there does {\it not} exist a surjective algebra homomorphism from
the affine Hecke algebra of type $G_2$ to the corresponding
Iwahori-Hecke subalgebra of type $G_2$.

\bigskip\noindent
{\bf Example 4. Type $G_2$.}
If $R$ is the root system of type $G_2$ then
$P=Q$ and $\Omega = \{1\}$.

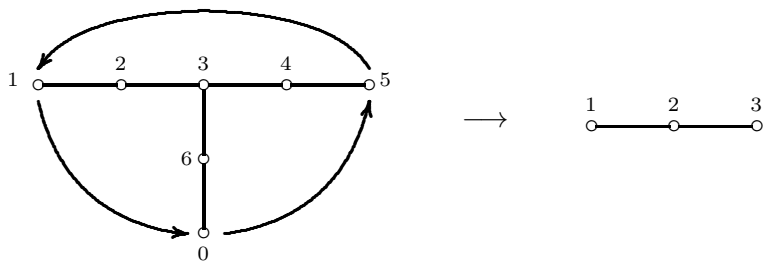
$$
%********************************************************************
% Dynkin diagram of type extended G_2
%********************************************************************
\beginpicture
\setcoordinatesystem units <1cm,1cm>       % sets scale
\setplotarea x from 2.8 to 4.2, y from -0.2 to 0.2 % sets plot size up
{\scriptsize
\multiput {$\circ$} at 3   0 *2 1 0 /      %  puts nodes in
\put {$^0$}     at 5 0.2   %
\put {$^2$}     at 4 0.2   %
\put {$^1$}     at 3 0.2   % label nodes with roots below
\linethickness=0.75pt                          % sets line thickness
\putrule from 3.03 0.055 to 3.97 0.055       % puts solid lines between 
nodes
\putrule from 3.03 -0.055 to 3.97 -0.055       %
\putrule from 3.05 0 to 3.95 0              %
\putrule from 4.05 0 to 4.95 0              %
}
\endpicture
\qquad
\beginpicture
\setcoordinatesystem units <1cm,1cm>       % sets scale
\setplotarea x from -1 to 1, y from -0.2 to 0.2 % sets plot size up
\put{$\longrightarrow$}    at 0 0
\put{$\times$}  at 0 0.01
\endpicture
\qquad
%********************************************************************
% Dynkin diagram of type G_2
%********************************************************************
\beginpicture
\setcoordinatesystem units <1cm,1cm>       % sets scale
\setplotarea x from 3.5 to 4.2, y from -0.2 to 0.2 % sets plot size up
{\scriptsize
\multiput {$\circ$} at 3   0 *1 1 0 /      % puts nodes in
\put {$^1$}     at 3 0.2   %
\put {$^2$}     at 4 0.2   % label nodes with roots above
\linethickness=0.75pt                      % sets line thickness
\putrule from 3.03 0.055 to 3.97 0.055  %
\putrule from 3.03 -0.055 to 3.97 -0.055  %
\putrule from 3.05 0 to 3.95 0         %   puts solid lines between nodes
}
\endpicture
$$

\prop
Let $\tilde H(G_2)$ be the affine Hecke algebra of
type $G_2$ as given by (4.6) and (4.10) and let $H(G_2)$ be the 
Iwahori-Hecke
subalgebra of type $G_2$ generated by $T_1$ and $T_2$.
There does not exist an algebra homomorphism
$\Phi\colon \tilde H(G_2)\to H(G_2)$ such that
$\Phi(T_i)=T_i$, for $1\le i\le 2$.
\pf
There is an irreducible representation
of $\tilde H(G_2)$ given by
$$\rho(T_1) =
\pmatrix{q &0\cr 0 &-q^{-1} \cr},
\qquad\hbox{and}\qquad
\rho(T_2) = {1\over q+q^{-1}}\pmatrix{
2-q^{-2} &q^2-1+q^{-2}\cr 3 &q^2-2\cr},$$
(see [Ra1, Theorem 6.11]).
We show that there does not exist a $2\times 2$ matrix $N$
which satisfies
$$N^2 = (q-q^{-1})N+1,\quad N\rho(T_2)N=\rho(T_2)N\rho(T_2)
\quad\hbox{and}\quad
N\rho(T_1)=\rho(T_1)N.$$
If $N$ exists then $N$ must be diagonal since
$N$ commutes with $\rho(T_1)$ and $\rho(T_1)$ is a diagonal
matrix with distinct eigenvalues.
The first equation shows that $N$ is invertible and
the second equation shows that $N$ is conjugate to
$\rho(T_2)$. It follows that $N$ must have one eigenvalue $q$
and one eigenvalue $-q^{-1}$.  Thus, either
$$N=\pmatrix{q &0\cr 0 &-q^{-1}\cr}
\qquad\hbox{or}\qquad
N=\pmatrix{-q^{-1} &0\cr 0 &q\cr}.$$
However, neither of these matrices satisfies the relation $N\rho(T_2)N
=\rho(T_2)N\rho(T_2).$  This  contradiction
shows that the representation $\rho$ cannot be extended to be a 
representation
of~$\tilde H(G_2)$.
\endpf

In spite of the fact, demonstrated by the previous example,
that there does not always exist a surjective algebra homomorphism
from the affine Hecke algebra onto its Iwahori-Hecke subalgebra,
there {\it are} interesting surjective homomorphisms from
affine Hecke algebras of exceptional type.

\bigskip\noindent
{\bf Example 5. Type $E_6$.}
For the root system of type $E_6$, $P/Q\cong \ZZ/3\ZZ$.
Let $\Omega = \{1,g,g^2\}$ where $g$ is as given by
(4.3) for the minuscule weight $\omega_1$ (see [Bou, p. 261]).
There are surjective algebra homomorphisms
$$
\beginpicture
\setcoordinatesystem units <1cm,1cm>                % sets scale
%********************************************************************
% Dynkin diagram of type E_6
%********************************************************************
\setplotarea x from 2.8 to 7.2, y from -1.7 to 0.4  % sets plot size
{\scriptsize
\multiput {$\circ$} at 7   0.5 *4 -1 0 /      % puts nodes in
\put{$\circ$} at 5 -0.4     %
\put{$\circ$} at 5 -1.3     %
\put {$^5$}   at 7 0.7   %
\put {$^4$}   at 6 0.7   %
\put {$^3$}   at 5 0.7   % label nodes
\put {$^2$}   at 4 0.7   %
\put {$^1$}   at 3 0.7   %
\put {$^6$}   at 4.8 -0.45   %
\put {$^0$}   at 4.8 -1.35   %
\setplotsymbol ({\rm .})
\linethickness=0.75pt                           % sets line thickness
\putrule from 6.05 0.5 to 6.95 0.5       %
\putrule from 5.05 0.5 to 5.95 0.5       %
\putrule from 4.05 0.5 to 4.95 0.5       %  puts solid lines between nodes
\putrule from 3.05 0.5 to 3.95 0.5       %
\putrule from 5 -.35 to 5 .45       %
\putrule from 5 -1.25 to 5 -.45     %
}
% end of dynkin diagram of type E_6
\endpicture
\qquad
\longrightarrow
\qquad
\beginpicture
\setcoordinatesystem units <1cm,1cm>                % sets scale
%********************************************************************
% Dynkin diagram of type A_5
%********************************************************************
\setplotarea x from 2.8 to 8.2, y from -0.8 to 0.8  % sets plot size
{\scriptsize
\multiput {$\circ$} at 7   0 *4 -1 0 /      % puts nodes in
\put {$^5$}   at 7 0.2   %
\put {$^4$}   at 6 0.2   %
\put {$^3$}   at 5 0.2   % label nodes
\put {$^2$}   at 4 0.2   %
\put {$^1$}   at 3 0.2   %
\setplotsymbol ({\rm .})
\linethickness=0.75pt                           % sets line thickness
\putrule from 6.05 0 to 6.95 0       %
\putrule from 5.05 0 to 5.95 0       %
\putrule from 4.05 0 to 4.95 0       %  puts solid lines between nodes
\putrule from 3.05 0 to 3.95 0       %
}
% end of dynkin diagram of type A_5
\endpicture
$$
$$\hskip.4in\matrix{
\Phi\colon &\tilde \cB_Q(E_6) &\longrightarrow &\cB(A_5) \cr
&T_0 &\longmapsto &T_5,\cr
&T_6 &\longmapsto &T_4,\cr
&T_i&\longmapsto &T_i, &\qquad &1\le i\le 5, \cr
}$$
and
\medskip
$$
\beginpicture
\setcoordinatesystem units <1cm,1cm>                % sets scale
%********************************************************************
% Dynkin diagram of type E_6
%********************************************************************
\setplotarea x from 2.8 to 7.2, y from -1.7 to 0.4  % sets plot size
{\scriptsize
\multiput {$\circ$} at 7   0.5 *4 -1 0 /      % puts nodes in
\put{$\circ$} at 5 -0.4     %
\put{$\circ$} at 5 -1.3     %
\put {$^5$}   at 7.2 0.5   %
\put {$^4$}   at 6 0.7   %
\put {$^3$}   at 5 0.7   % label nodes
\put {$^2$}   at 4 0.7   %
\put {$^1$}   at 2.7 0.5   %
\put {$^6$}   at 4.8 -0.45   %
\put {$^0$}   at 5 -1.6   %
\setplotsymbol ({\rm .})
\linethickness=0.75pt                           % sets line thickness
\putrule from 6.05 0.5 to 6.95 0.5       %
\putrule from 5.05 0.5 to 5.95 0.5       %
\putrule from 4.05 0.5 to 4.95 0.5       %  puts solid lines between nodes
\putrule from 3.05 0.5 to 3.95 0.5       %
\putrule from 5 -.35 to 5 .45       %
\putrule from 5 -1.25 to 5 -.45     %
\setquadratic \plot 7 0.3                      % draws automorphism
                     6.4 -.8                        %
                     5.25 -1.3 /                    %
\arrow <5pt> [.2,.67] from 6.98 0.2  to 7 0.3   %
\setquadratic \plot 3 0.3                      % draws automorphism
                     3.6 -.8                        %
                     4.75 -1.3 /                    %
\arrow <5pt> [.2,.67] from 4.7 -1.29 to 4.8 -1.3   %
\setquadratic \plot 7 0.7                      % draws automorphism
                     6.3 1.2                      %
                     5 1.4  /                   %
\setquadratic \plot 5 1.4                      %
                     3.7 1.2                      %
                     3 0.7 /                    %
\arrow <5pt> [.2,.67] from 3.08 0.8  to 3 0.7   %
}
% end of dynkin diagram of type E_6
\endpicture
\qquad
\longrightarrow
\qquad
\beginpicture
\setcoordinatesystem units <1cm,1cm>                % sets scale
%********************************************************************
% Dynkin diagram of type A_3
%********************************************************************
\setplotarea x from 2.8 to 6.2, y from -0.8 to 0.8  % sets plot size
{\scriptsize
\multiput {$\circ$} at 5   0 *2 -1 0 /      % puts nodes in
\put {$^3$}   at 5 0.2   % label nodes
\put {$^2$}   at 4 0.2   %
\put {$^1$}   at 3 0.2   %
\setplotsymbol ({\rm .})
\linethickness=0.75pt                           % sets line thickness
\putrule from 4.05 0 to 4.95 0       %  puts solid lines between nodes
\putrule from 3.05 0 to 3.95 0       %
}
% end of dynkin diagram of type A_3
\endpicture
$$
$$\hskip1.2in\matrix{
\Phi'\colon &\tilde \cB_P(E_6) &\longrightarrow &\cB(A_3) \cr
&g &\longmapsto &1, \cr
&T_0 &\longmapsto &T_1,\cr
&T_6 &\longmapsto &T_2,\cr
&T_5 &\longmapsto &T_1,\cr
&T_4 &\longmapsto &T_2,\cr
&T_i&\longmapsto &T_i, &\qquad &1\le i\le 3. \cr
}$$

\subsection{Example 6. Type $F_4$.}

For the root system $R$ of type $F_4$
we have $P=Q$ and $\Omega= \{1\}$.  There is a surjective homomorphism
$$
%********************************************************************
% Dynkin diagram of type extended F_4
%********************************************************************
\beginpicture
\setcoordinatesystem units <1cm,1cm>       % sets scale
\setplotarea x from 2.8 to 7.2, y from -0.2 to 0.2 % sets plot size up
{\scriptsize
\multiput {$\circ$} at 3   0 *4 1 0 /      %  puts nodes in
\put {$^0$}     at 7 0.2   %
\put {$^4$}     at 6 0.2   %
\put {$^3$}     at 5 0.2   %
\put {$^2$}     at 4 0.2   %
\put {$^1$}     at 3 0.2   % label nodes with roots below
\linethickness=0.75pt                          % sets line thickness
\putrule from 4.03 0.045 to 4.97 0.045       % puts solid lines between 
nodes
\putrule from 4.03 -0.045 to 4.97 -0.045       %
\putrule from 3.05 0 to 3.95 0              %
\putrule from 5.05 0 to 5.95 0              %
\putrule from 6.05 0 to 6.95 0              %
}
\endpicture
\qquad
\longrightarrow
\qquad
%********************************************************************
% Dynkin diagram of type A_2 x A_3
%********************************************************************
\beginpicture
\setcoordinatesystem units <1cm,1cm>       % sets scale
\setplotarea x from 2.8 to 7.2, y from -0.2 to 0.2 % sets plot size up
{\scriptsize
\multiput {$\circ$} at 3   0 *4 1 0 /      % puts nodes in
\put {$^0$}     at 7 0.2   %
\put {$^4$}     at 6 0.2   %
\put {$^3$}     at 5 0.2   %
\put {$^2$}     at 4 0.2   %
\put {$^1$}     at 3 0.2   % label nodes with roots below
\linethickness=0.75pt                      % sets line thickness
\putrule from 3.05 0 to 3.95 0              %
\putrule from 5.05 0 to 5.95 0              %
\putrule from 6.05 0 to 6.95 0              %
}
\endpicture
$$
$$\hskip1.7in\matrix{
\Phi\colon &\tilde \cB(F_4) &\longrightarrow &\cB(A_2\times A_3) \cr
&T_i&\longmapsto &T_i\,, &\qquad &0\le i\le 4. &\qquad\hbox{\qed}\cr
}$$

\bigskip
\sectno=6
\resultno=0
\centerline{\smallcaps Appendix: Clifford theory}~\medbreak

Let $R$ be an algebra over $\CC$ and let $G$ be a finite group acting by
automorphisms on $R$.  The {\it skew group ring} is
$$R\rtimes G
=\left\{ \sum_{g\in G} r_g g\ \Big|\ r_g\in R\right\}$$
with multiplication given by the distributive law and the relation
$$gr = g(r)g, \qquad\hbox{for $g\in G$ and $r\in R$.}$$
Let $N$ be a (finite dimensional) left $R$-module.  For each $g\in G$ 
define
an $R$-module $^gN$, which has the same underlying vector space $N$ but
such that
$$\hbox{$^gN$ has $R$-action given by}\qquad
r\circ n = g^{-1}(r)n,
\Aformula$$
for $r\in R$, $n\in N$.
If $W$ is an $R$-submodule of $N$ then
${}^gW$ is an $R$-submodule of ${}^gN$ and so ${}^gN$ is simple if and only
if $N$ is simple.  Thus there is an action of $G$ on the set of
simple $R$-modules.

Let $R^\lambda$ be a simple $R$-module.
The inertia group of $R^\lambda$ is
$$H = \{ h\in G\ |\ R^\lambda\cong\ ^hR^\lambda\}.
\Aformula$$
If $h\in H$ then Schur's lemma implies that the isomorphism
$R^\lambda\ \cong\  ^hR^\lambda$ is unique up to
constant multiples (since both $R^\lambda$ and $^hR^\lambda$ are simple).  
For
each $h\in H$ fix an isomorphism $\phi_h\colon R^\lambda\to\ ^{h^
{-1}}R^\lambda$.
Then,
as operators on $R^\lambda$,
$$\phi_h r=h(r) \phi_h,\quad\hbox{and}\quad \phi_g\phi_h =
\alpha(g,h)\phi_{gh},
\Aformula$$
where $\alpha(g,h)\in \CC^*$ are
determined by the choice of the isomorphisms
$\phi_h$.  The resulting function
$\alpha\colon H\times H \to \CC^*$ is called a {\it factor set} [CR, 8.32]
.

Let $(\CC H)_{\alpha^{-1}}$ be the algebra with basis $\{c_h\ |\ h\in H\}$
and multiplication given by
$$c_gc_h = \alpha(g,h)^{-1}c_{gh}, \qquad
\hbox{for $g,h\in H$.}
\Aformula$$
Let $H^\mu$ be a simple
$(\CC H)_{\alpha^{-1}}$-module.  Then putting
$$rh(m\otimes n) = r\phi_h m\otimes c_h n,\qquad
\hbox{for $r\in R$, $h\in H$, $m\in R^\lambda$, $n\in H^\mu$,}
\Aformula$$
defines an action of $R\rtimes H$ on $R^\lambda\otimes H^\mu$.

\Athm  (Clifford theory)  Let $R^\lambda$ be a simple $R$-module and let 
$H$ be
the inertia group of $R^\lambda$.  Let $H^\mu$ be a simple
$(\CC H)_{\alpha^{-1}}$-module where
$\alpha\colon H\times H\to \CC^*$ is the factor set determined by a choice 
of
isomorphisms $\phi_h\colon R^\lambda \to\ ^hR^\lambda$.  Define an action 
of
$R\rtimes H$ on $R^\lambda\otimes H^\mu$ as in (A.5) and define
$$RG^{\lambda,\mu}
=\Ind_{R\rtimes H}^{R\rtimes G}(R^\lambda\otimes H^\mu)
=(R\rtimes G) \otimes_{R\rtimes H} (R^\lambda\otimes H^\mu).$$
Then
\item{(a)}  $RG^{\lambda,\mu}$ is a simple $R\rtimes G$-module.
\item{(b)}  Every simple $R\rtimes G$-module is obtained by this
construction.
\item{(c)}  If $RG^{\lambda,\mu}\cong RG^{\nu,\gamma}$ then
$R^\lambda$ and $R^\nu$  are in the same $G$-orbit of simple $R$-modules
and $H^\mu\cong H^\gamma$ as $(\CC H)_{\alpha^{-1}}$-modules.
\pf
The proof of this theorem is as in [Mac2] except that the consideration
of the factor set $\alpha\colon H\times H\to \CC^*$ is necessary to
correct an error there.  We thank P. Deligne for pointing
this out to us.  A sketch of the proof is as follows.

Let $M$ be a simple $R\rtimes G$-module and let $R^\lambda$ be a simple
$R$-submodule of $M$.  Then $gR^\lambda \cong \ ^g R^\lambda$ as 
$R$-modules
and $M=\sum_{g\in G} gR^\lambda$
since the right hand side is an $R\rtimes G$-submodule of $M$.  Then
$$M = \sum_{g_i \in G/H} g_i N = \Ind_{R\rtimes H}^{R\rtimes G}(N),
\qquad \hbox{where}\quad N = \sum_{h\in H} hR^\lambda,$$
and the first sum is over a set $\{g_i\}$ of coset representatives of the
cosets
$G/H$.

The $R$-module $N$ is semisimple and by [Bou2]
$$N\cong R^\lambda\otimes H^\mu,
\Aformula$$
where $H^\mu = \Hom_R(R^\lambda,N)$.
It can be checked that the vector space $H^\mu$ has
a $(\CC H)_{\alpha^{-1}}$-action given by
$$(c_h\psi)(m) = \alpha(h,h^{-1}) h\psi(\phi_{h^{-1}}(m)),
\quad\hbox{for $h\in H$, $\psi\in \Hom_R(R^\lambda,N)$,}$$
where $c_h$ is as in (A.4).
Then, with $R\rtimes H$-action on $R^\lambda\otimes H^\mu$ given by (A.5),
the isomorphism in (A.7) is an isomorphism of $R\rtimes H$-modules (see
[CR, Thm. (11.17) (ii)]).

If $P$ is an $(\CC H)_{\alpha^{-1}}$-submodule of $H^\mu$ then
$R^\lambda\otimes P$ is an $R\rtimes H$-submodule of $R^\lambda\otimes 
H^\mu$
and $\Ind_{R\rtimes H}^{R\rtimes G}(R^\lambda\otimes P)$ is an
$R\rtimes G$-submodule of $M$.  Thus $H^\mu$ must be a simple
$(\CC H)_{\alpha^{-1}}$-module.

This argument shows that every simple $R\rtimes G$-module is
of the form $RG^{\lambda,\mu}$.
The uniqueness follows as in [Mac2, App.].
\endpf

\Aremark  A different choice $\psi_h\colon R^\lambda \to \ ^hR^\lambda$
of the isomorphisms in (A.3) may yield a factor set
$\beta\colon H\times H\to \CC^*$ which is different from the factor
set $\alpha$.  However, the algebras $(\CC H)_{\beta^{-1}}$ and
$(\CC H)_{\alpha^{-1}}$ are always isomorphic (a diagonal change of
basis suffices).

\Alemma    Define
$R^G=\{ r\in R  \ |\ \hbox{$g(r)=r$ for all $g\in G$}\}$ and let
$e = (1/|G|)\sum_{g\in G} g \in R\rtimes G.$
\item{(a)}  The map
$$\matrix{
\theta \colon & R^G & \longrightarrow & e(R\rtimes G) e \cr
& s & \longmapsto & se\cr
}$$
is a ring isomorphism.
\item{(b)}
Left multiplication
by elements of $R$ and the
action of $G$ by automorphisms
make $R$ into a left $R\rtimes G$-module.
Right multiplication makes $R$ a right $R^G$-module.
The rings $R\rtimes G$ and $e(R\rtimes G)e$ act
on $(R\rtimes G)e$ by left and right multiplication, respectively.
The map
$$\matrix{
\psi : & R & \cong & (R\rtimes G)e  \cr
& r & \longmapsto & re
}$$
is an isomorphism of $(R\rtimes G,R^G)$-bimodules.
\pf
(a)  If $r\in R^G$ then
$$ere = {1\over |G|}\sum_{g\in G} g(r)ge
={1\over |G|}\sum_{g\in G} re = re.$$
Thus the map $\theta$ is well defined and if $r,s\in R^G$ then $rese = 
rse$,
so $\theta$ is a homomorphism.
If $re=se$ then $r=s$ since $R\rtimes G$ is a free $R$-module with basis
$G$.  Thus $\theta$ is injective.  If $\sum_{g\in G} r_g g$ is a general
element of $R\rtimes G$ then
$$
e\left(\sum_{g\in G} r_g g \right) e =
\sum_{g,h\in G} h(r_g)hge
=\left(\sum_{g,h\in G} h(r_g)\right)e,$$
and, for each $g\in G$,  $\sum_{h\in G} h(r_g) \in R^G$.
So $\theta$ is surjective.

The proof of (b) is straightforward.
\endpf

Let $(\CC H)_{\alpha}$ be the algebra with basis $\{b_h\ |\ h\in H\}$
and multiplication given by
$$b_gb_h = \alpha(g,h)b_{gh}, \qquad
\hbox{for $g,h\in H$,}$$
and let $(\CC H)_{\alpha^{-1}}$ be as in (A.4).
Let $M$ be a $(\CC H)_\alpha$-module.  The
dual of $M$ is the $(\CC H)_{\alpha^{-1}}$-module given by the
vector space $M^*=\Hom(M,\CC)$ with action
$$(c_h \psi)(m) = \alpha(h,h^{-1})^{-1}\psi(b_{h^{-1}}m),
\qquad\hbox{for $h\in H$, $\psi\in M^*$.}$$
This is a $(\CC H)_{\alpha^{-1}}$ action since,
for all $g,h\in H$, $\psi\in M^*$,
$$\eqalign{
(c_gc_h \psi)(m)
&=\alpha(h,h^{-1})^{-1} \alpha(g,g^{-1})^{-1}
\psi( b_{h^{-1}} b_{g^{-1}} m ) \cr
&=\alpha(h,h^{-1})^{-1} \alpha(g,g^{-1})^{-1}
\alpha(h^{-1},g^{-1}) \psi( b_{(gh)^{-1}} m) \cr
&=\alpha(h,h^{-1})^{-1} \alpha(g,g^{-1})^{-1}
\alpha(h^{-1},g^{-1}) \alpha(gh,h^{-1}g^{-1})(c_{gh}\psi)(m) \cr
&=\alpha(g,h)^{-1}(c_{gh}\psi)(m), \cr}$$
where the last equality follows from the associativity of the product
$b_gb_hb_{h^{-1}}b_{g^{-1}}$ in $(\CC H)_{\alpha}$.
If $\rho\colon (\CC H)_\alpha\to \End(M)$ is the representation 
corresponding
to $M$ then the representation
$\rho^*\colon (\CC H)_{\alpha^{-1}}\to \End(M^*)$ corresponding to
$M^*$ is
$$\rho^*(c_h) = \alpha(h,h^{-1})^{-1}\rho(b_{h^{-1}})^t
=(\rho(b_h)^{-1})^t.
\Aformula$$

If $M$ is a $(\CC H)_{\alpha}$-module and
$N$ is a $(\CC H)_{\alpha^{-1}}$-module then
$M\otimes N$ is an $\CC H$-module with action defined by
$$h(m\otimes n)= b_h m\otimes c_h n,
\qquad\hbox{for $h\in H$, $m\in M$ and $n\in N$.}
\Aformula$$
The following lemma is a version of Schur's lemma which will be used in the
proof of Theorem~A.13.

\Alemma  Suppose that $M$ and $N$ are simple
$(\CC H)_\alpha$-modules and let $N^*$ be the
$(\CC H)_{\alpha^{-1}}$-module which is the dual of $N$.
Let $e_H = (1/|H|)\sum_{h\in H} h$.
Then
$$\dim(e_H(M\otimes N^*)) = \cases{
1, &if $M\cong N$, \cr
0, &otherwise.\cr}$$
\pf
Identify $M\otimes N^*$ with $\Hom(N,M)$.  Then, by (A.10) and (A.11),
the action of $\CC H$ on $\Hom(N,M)$ is given by,
$$hA = \rho(b_h)A\rho(b_h)^{-1},  \qquad
\hbox{for $h\in H$ and $A\in \Hom(N,M)$,}
$$
where $\rho\colon (\CC H)_{\alpha}\to \End(M)$ is the representation
corresponding to $M$.
If $A\in \Hom(N,M)$ and $g\in H$ then
$$e_HA = ge_HA = g(e_HA) = \rho(b_g)(e_H A)\rho(b_g)^{-1},$$
and so $\rho(b_g)(e_HA)=(e_HA)\rho(b_g)$ for all $g\in H$.
Then, by Schur's lemma, $e_HA=0$ if $M\not\cong N$ and
$e_HA$ is a constant if $M=N$.
\endpf

\Athm
Let $R^\lambda$ be a simple $R$-module and let $H$ be the inertia group of
$R^\lambda$.  The ring $R^G$ acts on $R^\lambda$ (by restriction) and
$(\CC H)_{\alpha}$ acts on $R^\lambda$ (by the $R$-module isomorphisms
$\phi_h\colon R^\lambda\ \cong \ ^h R^\lambda$ of (A.3))
and these two actions commute. Thus there is a decomposition
$$R^\lambda \cong \bigoplus_{\nu\in \hat H_\alpha}
R^{\lambda,\nu}\otimes (H^\nu)^*,$$
where $\hat H_\alpha$ is an index set for the simple $(\CC 
H)_\alpha$-modules,
$(H^\nu)^*$ is the dual of the simple $(\CC H)_{\alpha^{-1}}$-module 
$H^\nu$,
and $R^{\lambda,\nu}$ is an $R^G$-module.
\item{(a)}  If $R^{\lambda,\mu}\ne 0$ then it is a simple $R^G$-module.
\item{(b)}  Every simple $R^G$-module is isomorphic to some 
$R^{\lambda,\mu}$.
\item{(c)}  The nonzero $R^{\lambda,\mu}$ are pairwise nonisomorphic.
\pf
The setup of Lemma A.9(b) puts us in the situation of [Gr, \S 6.2].  If $e$
is the
idempotent used in Lemma A.9 then the
functor
$$\matrix{
\hbox{$R\rtimes G$-modules} &\longrightarrow &\hbox{$R^G$-modules}\cr
M &\longmapsto &eM\cr}$$
is an exact functor such that if $M$ is a simple $R\rtimes G$-module then
$eM$ is either $0$ or a simple $R^G$-module.  Furthermore, every simple
$R^G$-module arises as $eM$ for some simple $R\rtimes G$-module~$M$.

Let $RG^{\lambda,\mu}$ be a simple $R\rtimes G$-module as given by Theorem
A.6.  From the definition of
$RG^{\lambda,\mu}$ we obtain
$$\eqalign{
eRG^{\lambda,\mu}
&= e(R\rtimes G)\otimes_{R\rtimes H} (R^\lambda\otimes H^\mu)
\cr
&= e\otimes (R^\lambda\otimes H^\mu) = ee_H\otimes (R^\lambda\otimes H^\mu)
  \cr
&= e\otimes e_H(R^\lambda\otimes H^\mu),\cr}$$
where $e_H = (1/|H|)\sum_{h\in H} h$.  Using the decomposition in
the statement of the Theorem, we conclude that, as $R^G$-modules,
$$\eqalign{
eRG^{\lambda,\mu}
&= e\otimes e_H\left(\bigoplus_{\nu\in \hat H}
R^{\lambda,\nu}\otimes (H^\nu)^*\otimes H^\mu\right) \cr
&= e\otimes e_H\left(\bigoplus_{\nu\in \hat H}
R^{\lambda,\nu}\otimes e_H((H^\nu)^*\otimes H^\mu)\right) \cr
&\cong R^{\lambda,\mu}. \cr
}$$
The last isomorphism is a consequence of Lemma A.12.
The statement of the Theorem now follows from the results of
J.A. Green quoted above.

\Aremark  It follows from Theorem A.13 that $R^\lambda$ is
semisimple as an $R^G$-module and the action of $(\CC H)_\alpha$ on 
$R^\lambda$ generates
$\End_{R^G}(R^\lambda)$.

\vfill\eject
\centerline{\smallcaps References}

\medskip
\item{[Ar1]} {\smallcaps S.\ Ariki}
{\sl On the semisimplicity of the Hecke algebra of $(\ZZ/r\ZZ)\wr S_n$},
J.\ Algebra {\bf 169} (1994), 216--225.

\medskip
\item{[Ar2]} {\smallcaps S.\ Ariki}
{\sl Representation theory of a Hecke algebra of $G(r,p,n)$},  J.\ Algebra
{\bf 177} (1995), 164--185.

\medskip
\item{[Ar3]} {\smallcaps  S.\ Ariki},
{\sl On the decomposition numbers of the Hecke algebra of $G(m,1,n)$},
J.\  Math.\  Kyoto Univ.\ {\bf 36} (1996), 789--808.

\medskip
\item{[AK]} {\smallcaps S.\ Ariki and K.\ Koike}, {\it A Hecke algebra of
$(\ZZ/r\ZZ)\wr S_n$ and construction of its irreducible representations},
Adv.\ in Math.\ {\bf 106} (1994), 216--243.

\medskip
\item{[AM]} {\smallcaps S.\ Ariki and A.\ Mathas},
{\it On the number of simple modules of the Hecke algebras of 
type $G(r,p, n)$}, Math. Zeitschrift {\bf 233} (2000), 
no.\ 3, 601--623.

\medskip
\item{[Bou1]} {\smallcaps N.\ Bourbaki},
{\sl Groupes et alg\`ebres de Lie, Chapitres 4,5 et 6},
Elements de Math\'e\-matique, Hermann, Paris 1968.

\medskip
\item{[Bou2]} {\smallcaps N.\ Bourbaki},
{\sl Alg\`ebre, Chapitre 8},
Elements de Math\'e\-matique, Hermann, Paris 1958.

\medskip
\item{[BMM]} {\smallcaps M.\ Brou\'e, G.\ Malle, and J.\ Michel}, 
{\sl Repr\'esentations unipotents g\'en\'eriques et blocs des 
groupes r\'eductifs finis}, Ast\'erique {\bf 212}, 1993.

\medskip
\item{[Ch1]} {\smallcaps I.\ Cherednik},
{\it A new interpretation of Gel'fand-Tzetlin bases},
Duke Math. J. {\bf 54} (1987), 563--577.

\medskip
\item{[Ch2]} {\smallcaps I.\ Cherednik},
{\it Monodromy representations for generalized Knizhnik-Zamolodchikov
equations and Hecke algebras},
Publ.\ RIMS, Kyoto Univ.\ {\bf 27} (1991), 711--726.

\medskip
\item{[Cr]} {\smallcaps J.\ Crisp}, Ph.D. Thesis, University of Sydney, 
1997.

\medskip
\item{[CR]} {\smallcaps C.\ Curtis and I.\ Reiner}, {\sl
Methods of representation theory--with applications to finite groups and
orders} Vol. I, J.\ Wiley and Sons, 1981.

\medskip
\item{[DJ]} {\smallcaps R.\ Dipper and G.D.\ James},
{\it Blocks and idempotents of Hecke algebras of general linear groups},
Proc.\ London Math.\ Soc.\ (3) {\bf 54} (1987) 57--82.

\medskip
\item{[DJM]} {\smallcaps R.\ Dipper, G.D.\ James and G.E.\ Murphy},
{\it Hecke algebras of type $B_n$ at roots of unity},
Proc.\ London Math.\ Soc.\ (3) {\bf 70} (1995) 505--528.

\medskip
\item{[GL]} {\smallcaps M.\ Geck and S.\ Lambropoulou}, {\it Markov
traces and knot invariants related to Iwahori-Hecke algebras of type $B$}, 
J.  Reine Angew. Math. {\bf 482} (1997), 191--213.

\medskip
\item{[Gr]} {\smallcaps J.A. Green}, {\sl Polynomial representations of
$GL_n$}, Lecture Notes in Mathematics {\bf 830}, Springer-Verlag, New York,
1980.

\medskip
\item{[Gy]}  {\smallcaps A. Gyoja}, 
{\it A q-analogue of Young symmetrizer},
Osaka J.\ Math.\ {\bf 23} (1986), 841--852.

\medskip
\item{[HR]} {\smallcaps T.\ Halverson and A.\ Ram},
{\it Murnaghan-Nakayama rules for characters of Iwahori-Hecke algebras of 
the
complex reflection groups $G(r,p,n)$},
Canadian J.\ Math.\ {\bf 50} (1998), 167--192.

\medskip
\item{[Hf]} {\smallcaps P.N.\ Hoefsmit},
{\it Representations of Hecke algebras of finite groups with $BN$-pairs
of classical type},
Ph.D.\ Thesis, University of British Columbia, 1974.

\medskip
\item{[IM]} {\smallcaps N.\ Iwahori and H. Matsumoto},
{\it On some Bruhat decomposition and the structure of the 
Hecke rings of $p$-adic Chevalley groups}, 
Publ.\ Math.\ I.H.E.S. {\bf 25} (1965), 5--48.

\medskip
\item{[Jo]} {\smallcaps V.F.R.\ Jones}, {\it A quotient of the affine Hecke
algebra in the Brauer algebra}, Enseign. Math. (2) {\bf 40} (1994), 
313--344.

\medskip
\item{[KL]} {\smallcaps D.\ Kazhdan and G.\ Lusztig},
{\it Proof of the Deligne-Langlands conjecture for Hecke algebras},
Invent. Math. {\bf 87} (1987), 153--215.

\medskip
\item{[Lu]} {\smallcaps G.\ Lusztig},
{\it Affine Hecke algebras and their graded version}, J. Amer. Math.
Soc. {\bf 2} (1989), 599--635.

\medskip
\item{[Mac]} {\smallcaps I.G.\ Macdonald},
{\sl Symmetric functions and Hall polynomials}, Second edition, Oxford
University Press, New York, 1995.

\medskip
\item{[Mac2]} {\smallcaps  I.G.\ Macdonald}, {\it Polynomial functors and
wreath products}, J.\ Pure and Appl.\ Algebra {\bf 18} (1980) 173--204.

\medskip
\item{[Mac3]} {\smallcaps I.G.\ Macdonald},
{\it Affine Hecke algebras and orthogonal polynomials},
S\'eminaire Bourbaki, 47\`eme ann\'ee, ${\rm n}^{\rm o}$ 797, 1994--95,
Ast\'erisque {\bf 237} (1996), 189--207.
\medskip

\medskip
\item{[Ra1]} {\smallcaps A.\ Ram},
{\it Seminormal representations of Weyl groups and Iwahori-Hecke algebras},
Proc.\ London Math.\ Soc.\ (3) {\bf 75} (1997), 99-133.

\medskip
\item{[Ra2]} {\smallcaps A.\ Ram},
{\it Calibrated representations of affine Hecke algebras}, preprint 1998,
results to appear in J.\ Algebra under the title {\it Affine Hecke algebras
and generalized standard Young tableaux}.

\medskip
\item{[Ra3]} {\smallcaps A. \ Ram},
{\it Standard Young tableaux for finite root systems}, preprint 1998,
results to appear in J.\ Algebra under the title {\it Affine Hecke algebras
and generalized standard Young tableaux}.

\medskip
\item{[Ra4]} {\smallcaps A.\ Ram},
{\it Irreducible representations of rank two affine Hecke algebras},
preprint 1998.

\medskip
\item{[Ra5]} {\smallcaps A.\ Ram}, {\it Skew shape representations
are irreducible}, preprint 1998.

\medskip
\item{[Re]} {\smallcaps M. Reeder}, {\it Isogenies of Hecke algebras
and a Langlands correspondence for ramified principal series 
representations}, Represent.\ Theory {\bf 6} (2002), 101--126.

\medskip
\item{[RR]} {\smallcaps A.\ Ram and J.\ Ramagge},
{\sl Langlands parameters for affine Hecke algebra representations},
in preparation.

\medskip
\item{[Sp]} {\smallcaps W. Specht}, {\sl Eine Verallgemeinerung der
symmetrischen Gruppe}, Schriften Math. Seminar Berlin {\bf 1} (1932), 
1--32.

\medskip
\item{[ST]} {\smallcaps G. C.\ Shephard and J.A.\ Todd}, 
{\sl Finite unitary reflection groups}, Canadian Journal of Math.\ {\bf 6} 
(1954), 274--304.

\medskip
\item{[Wz]}  {\smallcaps H. Wenzl},
{\it Hecke algebras of type $A\sb n$ and subfactors},
Invent. Math. {\bf 92} (1988), 349--383.

\medskip
\item{[Yg1]}  {\smallcaps A.\ Young},
{\it On quantitative substitutional analysis} (fifth paper),
Proc. London Math. Soc. (2) {\bf 31} (1929), 273--288.

\medskip
\item{[Yg2]}  {\smallcaps A.\ Young},
{\it On quantitative substitutional analysis} (sixth paper),
Proc. London Math. Soc. (2) {\bf 34} (1931), 196--230.

\vfill\eject
\end